\documentclass[leqno, 12pt]{article}
\usepackage{amsmath,amsfonts,amsthm,amssymb,indentfirst}

\newtheorem{theorem}{Theorem}
\newtheorem{lemma}[theorem]{Lemma}
\newtheorem{definition}[theorem]{Definition}
\newtheorem{corollary}[theorem]{Corollary}
\newtheorem{proposition}[theorem]{Proposition}

\newcommand{\E}{\mathrm{Self}}

\newcommand{\Sym}{\mathrm{Sym}}

\newcommand{\preq}{\preccurlyeq}

\makeatletter
\newcommand{\xlabel}{\stepcounter{equation}
  \gdef\@currentlabel{\p@equation\theequation}{\rm(\@currentlabel)}}
\makeatother
\newenvironment{xlist}
  {\begin{list}{\xlabel}
    {\setlength{\rightmargin}{20pt}
     \setlength{\leftmargin}{37pt}
     \setlength{\labelsep}{20pt}
     \setlength{\labelwidth}{20pt}}}
  {\end{list}}

\begin{document}

\title{Generating self-map monoids of infinite sets\thanks{
(secondary).}}
\author{Zachary Mesyan}

\maketitle

\begin{abstract}
Let $\Omega$ be a countably infinite set, $S = \Sym(\Omega)$ the
group of permutations of $\Omega$, and $E = \E(\Omega)$ the monoid
of self-maps of $\Omega$. Given two subgroups $G_1, G_2 \subseteq
S$, let us write $G_1 \approx_S G_2$ if there exists a finite
subset $U \subseteq S$ such that the groups generated by $G_1 \cup
U$ and $G_2 \cup U$ are equal. Bergman and Shelah showed that the
subgroups which are closed in the function topology on $S$ fall
into exactly four equivalence classes with respect to $\approx_S$.
Letting $\approx$ denote the obvious analog of $\approx_S$ for
submonoids of $E$, we prove an analogous result for a certain
class of submonoids of $E$, from which the theorem for groups can
be recovered. Along the way, we show that given two subgroups
$G_1, G_2 \subseteq S$ which are closed in the function topology
on $S$, we have $G_1 \approx_S G_2$ if and only if $G_1 \approx
G_2$ (as submonoids of $E$), and that $\mathrm{cl}_S (G) \approx
\mathrm{cl}_E (G)$ for every subgroup $G \subseteq S$ (where
$\mathrm{cl}_S (G)$ denotes the closure of $G$ in the function
topology in $S$ and $\mathrm{cl}_E (G)$ its closure in the
function topology in $E$).
\end{abstract}

\section{Introduction}

Let $\Omega$ be a countably infinite set, let $S = \Sym(\Omega)$
denote the group of all permutations of $\Omega$, and let $E =
\E(\Omega)$ denote the monoid of self-maps of $\Omega$. (Here
``$E$" stands for ``endomap.") Given two subgroups $G_1, G_2
\subseteq S$, let us write $G_1 \approx_S G_2$ if there exists a
finite subset $U \subseteq S$ such that the group generated by
$G_1 \cup U$ is equal to the group generated by $G_2 \cup U$.
In~\cite{B&S} Bergman and Shelah show that the subgroups of $S$
that are closed in the function topology on $S$ fall into exactly
four equivalence classes with respect to the above equivalence
relation. (A subbasis of open sets in the function topology on $S$
is given by the sets $\{f \in S : (\alpha)f = \beta\}$
$(\alpha,\beta \in \Omega).$ This topology is discussed in more
detail in Sections~\ref{topology section} and~\ref{groups
section}.) In this note we investigate properties of an analogous
equivalence relation $\approx$ defined for monoids. Two submonoids
$M_1, M_2 \subseteq E$ will be considered equivalent if and only
if there exists a finite set $U \subseteq E$ such that the monoid
generated by $M_1 \cup U$ is equal to the monoid generated by $M_2
\cup U$.

We will show that given two subgroups $G_1, G_2 \subseteq S$ that
are closed in the function topology on $S$, we have $G_1 \approx_S
G_2$ if and only if $G_1 \approx G_2$ (as submonoids of $E$).
Writing $\mathrm{cl}_S (G)$ for the closure of the subgroup $G
\subseteq S$ in the function topology in $S$ and $\mathrm{cl}_E
(G)$ for its closure in the function topology in $E$, we will show
that $\mathrm{cl}_S (G) \approx \mathrm{cl}_E (G)$.

Our main goal will be to classify into equivalence classes certain
closed submonoids of $E$. In general, unlike the case of groups,
the closed submonoids of $E$ fall into infinitely many equivalence
classes. To see this, let $\Omega = \omega$, the set of natural
numbers, and for each positive integer $n$ let $M_n$ be the
submonoid of $E$ generated by all maps whose images are contained
in $\{0, 1, \dots, n\}$. It is easy to see that such submonoids
are closed in the function topology. Now, let $n$ be a positive
integer and consider a monoid word $g_1 g_2 \dots g_k$ in elements
of $E$, where at least one of the $g_i \in M_n$. Then the image of
$\Omega$ under $g_1 g_2 \dots g_k$ has cardinality at most $n$.
Hence if for some finite subset $U \subseteq E$ we have that
$M_{n+1} \subseteq \langle M_n \cup U \rangle$ (the monoid
generated by $M_n \cup U$), then all the elements of $M_{n+1}$
whose images have cardinality $n+1$ must be in $\langle U
\rangle$. This is impossible, since there are uncountably many
such elements. Therefore, if $n \neq m$ are two positive integers,
then $M_n \not\approx M_m$.

Classifying all closed submonoids of $E$ into equivalence classes
appears to be a very difficult task, so we will for the most part
focus on submonoids that have large stabilizers (i.e., submonoids
$M$ such that for any finite set $\Sigma \subseteq \Omega$ the
pointwise stabilizer of $\Sigma$ in $M$ is $\approx M$) and the
property that the subset of a.e.\ injective maps is dense (with
respect to the function topology). Examples of submonoids that
have large stabilizers and dense sets of a.e.\ injective maps
include subgroups of $S \subseteq E$, as well as their closures in
the function topology in $E$. Another class of examples of monoids
with these properties arises from preorders on $\Omega$. Given
such a preorder $\rho$, let $E(\rho)$ denote the submonoid of $E$
consisting of all maps $f$ such that for all $\alpha \in \Omega$,
one has $(\alpha, (\alpha)f) \in \rho$. We will show that the
submonoids $E(\rho)$ have large stabilizers and dense sets of
a.e.\ injective maps. (They are also closed in the function
topology on $E$.)

The Bergman-Shelah theorem can be recovered from our
classification of the submonoids described above. A surprising
aspect of this classification is that the submonoids in question
fall into five equivalence classes, rather than the four predicted
by the Bergman-Shelah theorem. Throughout, we will also give a
number of examples demonstrating various unusual features of our
equivalence relation.

Conditions under which a submonoid $M \subseteq E$ satisfies $M
\approx E$ have been considered before. For instance, Howie,
Ru\v{s}kuc, and Higgins show in~\cite{HRH} that $S \approx E$. (We
give a shorter proof of this fact below.) These authors and
Mitchell also exhibit various submonoids that are $\not\approx E$
in~\cite{HHMR}. Related questions, but with other kinds of objects
in place of $E$, are discussed in~\cite{CMM}, \cite{HHMR},
and~\cite{ZM2}, as well as in papers referenced therein.

\section{Chains}

Let $\Omega$ be an arbitrary infinite set, and set $E =
\E(\Omega)$, the monoid of self-maps of $\Omega$. Elements of $E$
will be written to the right of their arguments. If $U \subseteq
E$ is a subset, then we will write $\langle U \rangle$ to denote
the submonoid generated by $U$. The cardinality of a set $\Gamma$
will be denoted by $|\Gamma|$. If $\Sigma \subseteq \Omega$ and $U
\subseteq E$ are subsets, let $U_{\{\Sigma\}} = \{f \in U :
(\Sigma) f \subseteq \Sigma \}$.

\begin{definition}
Let $M$ be a monoid that is not finitely generated. Then the {\em
cofinality} $c(M)$ of $M$ is the least cardinal $\kappa$ such that
$M$ can be expressed as the union of an increasing chain of
$\kappa$ proper submonoids.
\end{definition}

The main goal of this section is to show that $c(E) > |\Omega|$,
which will be needed later on. This section is modeled on Sections
1 and 2 of~\cite{GB}, where analogous statements are proved for
the group of all permutations of $\Omega$. Ring-theoretic analogs
of these statements are proved in~\cite{ZM}.

\begin{lemma}\label{g,h}
Let $U \subseteq E$ and $\Sigma \subseteq \Omega$ be such that
$|\Sigma| = |\Omega|$ and the set of self-maps of $\Sigma$ induced
by $U_{\{\Sigma\}}$ is all of $\E(\Sigma)$. Then $E = gUh$, for
some $g, h \in E$.
\end{lemma}

\begin{proof}
Let $g \in E$ be a map that takes $\Omega$ bijectively to
$\Sigma$, and let $h \in E$ be a map whose restriction to $\Sigma$
is the right inverse of $g$. Then $E = gU_{\{\Sigma\}}h$.
\end{proof}

We will say that $\Sigma \subseteq \Omega$ is a \emph{moiety} if
$|\Sigma| = |\Omega| = |\Omega \backslash \Sigma|$. A moiety
$\Sigma \subseteq \Omega$ is called \emph{full} with respect to $U
\subseteq E$ if the set of self-maps of $\Sigma$ induced by
members of $U_{\{\Sigma\}}$ is all of $\E(\Sigma)$. The following
two results are modeled on group-theoretic results in~\cite{M&N}.

\begin{lemma}[{\rm cf.\ \cite[Lemma~3]{GB}}]\label{diagonal}
Let $(U_i)_{i \in I}$ be any family of subsets of $E$ such that
$\bigcup_{i \in I} U_i = E$ and $|I| \leq |\Omega|$. Then $\Omega$
contains a full moiety with respect to some $U_i$.
\end{lemma}

\begin{proof}
Since $|\Omega|$ is infinite and $|I| \leq |\Omega|$, we can write
$\Omega$ as a union of disjoint moieties $\Sigma_i$, $i \in I$.
Suppose that there are no full moieties with respect to $U_i$ for
any $i \in I$.  Then in particular, $\Sigma_i$ is not full with
respect to $U_i$ for any $i \in I$.  Hence, for every $i \in I$
there exists a map $f_i \in \E(\Sigma_i)$ which is not the
restriction to $\Sigma_i$ of any member of $(U_i)_{\{\Sigma_i\}}$.
Now, if we take $f \in E$ to be the map whose restriction to each
$\Sigma_i$ is $f_i$, then $f$ is not in $U_i$ for any $i \in I$,
contradicting $\bigcup_{i \in I} U_i = E$.
\end{proof}

\begin{proposition}\label{chains}
$c(E) > |\Omega|$.
\end{proposition}

\begin{proof}
Suppose that $(M_i)_{i \in I}$ is a chain of submonoids of $E$
such that $\bigcup_{i \in I} M_i = E$ and $|I| \leq |\Omega|$. We
will show that $E = M_i$ for some $i \in I$.

By the preceding lemma, $\Omega$ contains a full moiety with
respect to some $M_i$. Thus, Lemma~\ref{g,h} implies that $E =
\langle M_i \cup \{g, h \} \rangle$ for some $g, h \in E$. But, by
the hypotheses on $(M_i)_{i \in I}$, $M_i \cup \{g, h \} \subseteq
M_j$ for some $j \in I$, and hence $E = \langle M_i \cup \{g, h \}
\rangle \subseteq M_j$, since $M_j$ is a submonoid.
\end{proof}

This result is proved by a very different method in~\cite{GB2}.

\section{Equivalence classes}

Throughout this note we will be primarily interested in submonoids
of $E = \E(\Omega)$. However, we begin this section with a
definition applicable to submonoids of an arbitrary monoid.

\begin{definition}
Let $M$ be a monoid, $\kappa$ an infinite cardinal, and $M_1,\
M_2$ submonoids of $M$.  We will write $M_1 \preq_{\kappa,M} M_2$
if there exists a subset $U \subseteq M$ of cardinality $<\kappa$
such that $M_1 \subseteq \langle M_2 \cup U \rangle$. If $M_1
\preq_{\kappa,M} M_2$ and $M_2 \preq_{\kappa,M} M_1,$ we will
write $M_1 \approx_{\kappa,M} M_2$, while if $M_1 \preq_{\kappa,M}
M_2$ and $M_2 \not\preq_{\kappa,M} M_1$, we will write $M_1
\prec_{\kappa,M} M_2$.  The subscripts $M$ and $\kappa$ will be
omitted when their values are clear from the context.
\end{definition}

It is clear that $\preq_{\kappa,M}$ is a preorder on submonoids of
$M$, and hence $\approx_{\kappa,M}$ is an equivalence relation.
This equivalence relation and many of the results below are
modeled on those in~\cite{B&S}. Ring-theoretic analogs of these
results can be found in~\cite{ZM2}.

We record the following result for future use.

\begin{theorem}[{\rm Sierpi\'{n}ski, cf.\ \cite{WS}}, \cite{Banach},
\cite{HHMR}]\label{2-gen} Every countable subset of $E$ is
contained in a subsemigroup generated by two elements of $E$.
\end{theorem}

\begin{proposition}\label{reductions}
Let $M_1, M_2 \subseteq E$ be submonoids.\\[2pt]
{\rm(i)} $M_1 \preq_{\aleph_0} M_2$ if and only if $M_1
\preq_{\aleph_1} M_2$ {\rm(}and hence $M_1 \approx_{\aleph_0} M_2$
if and only if $M_1 \approx_{\aleph_1} M_2)$.\\[2pt]
{\rm(ii)} $M \approx_{\aleph_0} E$ if and only if $M
\approx_{|\Omega|^+} E$ $($where $|\Omega|^+$ is the successor
cardinal of $|\Omega|$$)$.
\end{proposition}

\begin{proof}
(i) follows from Theorem~\ref{2-gen}.  (ii) follows from
Proposition~\ref{chains}. For, if $M \approx_{|\Omega|^+} E$, then
among subsets $U \subseteq E$ of cardinality $\leq |\Omega|$ such
that $\langle M \cup U \rangle = E$, we can choose one of least
cardinality. Let us write $U = \{f_i : i \in |U|\}$. Then the
submonoids $M_i = \langle M \cup \{f_j : j < i \} \rangle$ $(i \in
|U|)$ form a chain of $\leq |\Omega|$ proper submonoids of $E$. If
$|U|$ were infinite, this chain would have union $E$,
contradicting Proposition~\ref{chains}. Hence $U$ is finite, and
$M \approx_{\aleph_0} E$.
\end{proof}

\begin{definition}\label{leq}
Let $E_\leq \subseteq \E(\omega)$ denote the submonoid of all maps
decreasing with respect to the usual ordering of $\omega$.
Specifically, $f \in E_\leq$ if and only if for all $\alpha \in
\omega$, $(\alpha) f \leq \alpha$
\end{definition}

The following result will be our main tool for separating various
equivalence classes of submonoids of $\E(\Omega)$ throughout the
paper.

\begin{theorem}\label{fear}
\begin{enumerate}
\item[\rm{(i)}] Let $\kappa$ be a regular infinite cardinal
$\leq |\Omega|$ and $T \subseteq \E(\Omega)$ a subset. If
$|(\alpha)T| < \kappa$ for all $\alpha \in \Omega$, then $\langle
T \rangle \prec_{|\Omega|^+} \E(\Omega)$.

\item[\rm{(ii)}] Let $T \subseteq \E(\omega)$ be a subset. If there exists
a finite $\lambda$ such that $|(\alpha)T| \leq \lambda$ for all
$\alpha \in \omega$, then $\langle T \rangle \prec_{\aleph_0}
E_\leq$.
\end{enumerate}
\end{theorem}

\begin{proof}
(i) Let $U \subseteq \E(\Omega)$ be a subset of cardinality $\leq
|\Omega|$. We will show that $\E(\Omega) \not\subseteq \langle T
\cup U \rangle$. Without loss of generality we may assume that $1
\in T \cap U$. For all $j \in \omega$ and $u_0, \dots, u_j \in U$,
we define
\begin{xlist}\item
$B(u_0, \dots, u_j) = \{u_0t_0 \dots u_jt_j : t_0, \dots, t_j \in
T\}.$
\end{xlist}
Then the monoid $\langle T \cup U \rangle$ can be written as the
union of the sets $B(u_0, \dots, u_j)$ (using the assumption that
$1 \in T \cap U$).

Next, we show by induction on $j$ that for all $\alpha \in
\Omega$, $j \in \omega$, and $u_0, \dots, u_j \in U$, we have
$|(\alpha)B(u_0, \dots, u_j)| < \kappa$. If $j = 0$,
$|(\alpha)B(u_0)| = |((\alpha)u_0)T| < \kappa$, by our hypothesis
on $T$. Now, $(\alpha)B(u_0, \dots, u_{j+1}) =
((\Sigma)u_{j+1})T$, where $\Sigma = (\alpha)B(u_0, \dots, u_j)$.
Assuming inductively that $|\Sigma| < \kappa$, and hence
$|(\Sigma)u_{j+1}| < \kappa$, the set $(\alpha)B(u_0, \dots,
u_{j+1})$ can be written as the union of $< \kappa$ sets of
cardinality $< \kappa$. By the regularity of $\kappa$,
$|(\alpha)B(u_0, \dots, u_{j+1})| < \kappa$.

Let us write $\Omega = \bigcup_{j\in\omega} \Omega_j$, where the
union is disjoint and each $\Omega_j$ has cardinality $|\Omega|$.
Also, for each $n \in \omega$ let $h_n : \Omega_n \rightarrow
\prod_{j=0}^n U$ be a surjection. By the previous paragraph, there
is a map $f \in \E(\Omega)$ such that for all $\alpha \in
\Omega_j$, $(\alpha)f \in \Omega \setminus (\alpha)B(u_0, \dots,
u_j)$, where $(u_0, \dots, u_j) = (\alpha)h_j$. We conclude the
proof by showing that $f \notin \langle T \cup U \rangle$.

Suppose that $f \in \langle T \cup U \rangle$. Then $f \in B(u_0,
\dots, u_j)$ for some $j \in \omega$ and $u_0, \dots, u_j \in U$.
Let $\alpha \in \Omega_j$ be such that $(u_0, \dots, u_j) =
(\alpha)h_j$. Then $(\alpha)f \in (\alpha)B(u_0, \dots, u_j)$,
contradicting our definition of $f$. Hence $f \notin \langle T
\cup U \rangle$.

(ii) Let $T \subseteq \E(\omega)$ and $\lambda \in \omega$ be such
that $|(\alpha)T| \leq \lambda$ for all $\alpha \in \omega$, and
let $U \subseteq \E(\omega)$ be a finite subset. We note that for
all $\alpha \in \omega$, $|(\alpha)(T \cup U)| \leq \lambda + |U|
< \aleph_0$. Hence $T \cup U$ satisfies our hypotheses on $T$.
Therefore, to show that $E_\leq \not\subseteq \langle T \cup U
\rangle$ for all $T$ and $U$, it suffices to show that $E_\leq
\not\subseteq \langle T \rangle$ for all $T$ as in the statement.

Let $f \in E_\leq$ be any element such that
\begin{xlist}\item
$(\lambda^j+1)f \notin (\lambda^j+1)T^j$ for all $j \geq 1$,
\end{xlist}
where $T^j = \{t_0 \dots t_{j-1} : t_0, \dots, t_{j-1} \in T\}$.
Such a map exists, since $|(\lambda^j+1)T^j| \leq \lambda^j$,
while there are $\lambda^j+1$ possible values for
$(\lambda^j+1)f$. Then $f \notin T^j$ for all $j \geq 1$, and
hence $f \notin \langle T \rangle$. It remains to be shown that
$\langle T \rangle \preq_{\aleph_0} E_\leq$.

We can find a set $\omega' \subseteq \omega$ and a collection of
disjoint sets $\Delta_\alpha \subseteq \omega$ ($\alpha \in
\omega'$) such that $|\omega'| = |\omega|$, $|\Delta_\alpha| =
\lambda$, and $\Delta_\alpha \subseteq (\alpha)E_\leq$.
(Specifically, we can take $\omega' = \{\lambda, 2\lambda,
3\lambda, \dots \}$ and $\Delta_{i\lambda} = \{(i-1)\lambda,
(i-1)\lambda + 1, (i-1)\lambda + 2, \dots, i\lambda - 1\}$ for $i
\geq 1$.) Now, let $g \in \E(\omega)$ be an injective map from
$\omega$ to $\omega'$, and let $h \in \E(\omega)$ be a map that
takes each $\Delta_{(\alpha) g}$ ($\alpha \in \omega'$) onto
$(\alpha) T$. Then $T \subseteq g E_\leq h$, and hence $\langle T
\rangle \subseteq \langle E_\leq \cup \{g, h\} \rangle$.
\end{proof}

While we will not use the following two results in the future,
they are of interest in their own right.

\begin{corollary}\label{union fear}
Let $\kappa$ be a regular infinite cardinal $\leq |\Omega|$ and
$\{T_i\}_{i\in \kappa}$ subsets of $E = \E(\Omega)$ such that
$|(\alpha)T_i| < \kappa$ for all $\alpha \in \Omega$ and $i\in
\kappa$. Then $\langle \bigcup_{i\in \kappa} T_i \rangle
\prec_\kappa E$.
\end{corollary}

\begin{proof}
Suppose that $\langle \bigcup_{i\in \kappa} T_i \rangle
\approx_\kappa E$. Then there is a set $U \subseteq E$ of
cardinality $< \kappa$ such that $E = \langle \bigcup_{i\in
\kappa} T_i \cup U \rangle$. For each $j \in \kappa$ let $N_j =
\langle \bigcup_{i \leq j} T_i \cup U \rangle$. Then $\langle
\bigcup_{i\in \kappa} T_i \cup U \rangle = \bigcup_{i \in \kappa}
N_i$, and so $E = \bigcup_{i \in \kappa} N_i$. Hence, by
Proposition~\ref{chains}, $E = N_n$ for some $n \in \kappa$. For
each $\alpha \in \Omega$, set $|(\alpha)\bigcup_{i \leq n} T_i| =
\lambda_\alpha$. Then each $\lambda_\alpha < \kappa$, and for all
$\alpha \in \Omega$, we have that $|(\alpha)\bigcup_{i \leq n} T_i
\cup U| \leq \lambda_\alpha + |U| < \kappa$. Thus, $E = N_n$
contradicts Theorem~\ref{fear}.
\end{proof}

In view of Proposition~\ref{reductions}(ii), $\prec_\kappa$ can be
replaced with $\prec_{|\Omega|^+}$ in the above corollary. This
result can be viewed as a generalization to arbitrary $\Omega$ of
\cite[Corollary~2.2]{HHMR}. In~\cite{HHMR} a subset $T \subseteq
\E(\omega)$ is said to be {\em dominated} (by $U$) if there exists
a countable subset $U \subseteq \E(\omega)$ having the property
that for each $f \in T$ there exists $h \in U$ such that
$(\alpha)f \leq (\alpha)h$ for all $\alpha \in \omega$. Rewritten
using our notation, Corollary 2.2 states that if $T \subseteq
\E(\omega)$ is a dominated subset, then $\langle T \rangle
\prec_{\aleph_1} \E(\omega)$. This can be deduced from
Corollary~\ref{union fear} as follows. Let $U \subseteq
\E(\omega)$ be a countable subset that dominates $T \subseteq
\E(\omega)$, and write $U = \{h_i : i \in \omega \}$. Then $T =
\bigcup_{i \in \omega} T_i$, where $T_i \subseteq T$ is a subset
dominated by $\{h_i\}$. For all $\alpha \in \omega$, we have that
$|(\alpha)T_i| \leq (\alpha)h_i + 1 < \aleph_0$.
Corollary~\ref{union fear} then implies that $\langle T \rangle  =
\langle \bigcup_{i \in \omega} T_i \rangle \prec_{\aleph_0}
\E(\omega)$, which is equivalent to $\langle T \rangle
\prec_{\aleph_1} \E(\omega)$.

The next result shows that if the $\kappa$ in the statement of
Theorem~\ref{fear} is assumed to be uncountable, then a stronger
conclusion can be obtained (with less work).

\begin{proposition}
Let $\kappa$ be a regular uncountable cardinal $\leq |\Omega|$ and
$T \subseteq E$ a subset.
\begin{enumerate}
\item[\rm{(i)}]  If $|(\alpha)T| < \kappa$ for all $\alpha \in \Omega$,
then $|(\alpha)\langle T \rangle| < \kappa$ for all $\alpha \in
\Omega$.

\item[\rm{(ii)}] If there exists $\aleph_0 \leq \lambda < \kappa$ such that
$|(\alpha)T| \leq \lambda$ for all $\alpha \in \Omega$, then
$|(\alpha)\langle T \rangle| \leq \lambda$ for all $\alpha \in
\Omega$.
\end{enumerate}
\end{proposition}

\begin{proof}
For each $\beta \in \Omega$, let $\lambda_\beta < \kappa$ be such
that $|(\beta) T| \leq \lambda_\beta$. Also, let $\alpha \in
\Omega$ be any element. Then, by definition,
\begin{xlist}\item
$\displaystyle (\alpha)\langle T \rangle = \bigcup_{j=1}^\infty
(\alpha)T^j$.
\end{xlist}
We claim that $(\alpha)T^j < \kappa$ for all $j \geq 1$. This is
true, by hypothesis, for $j = 1$. Assuming inductively that
$\Sigma = (\alpha)T^{j-1}$ has cardinality $< \kappa$, we have
\begin{xlist}\item
$\displaystyle |(\alpha)T^j| = |(\Sigma)T| = |\bigcup_{\sigma \in
\Sigma} \{(\sigma)f : f \in T\}| \leq \sum_{\sigma \in \Sigma}
\lambda_\sigma$.
\end{xlist}
This sum has $< \kappa$ summands, each $< \kappa$; therefore
$|(\alpha)T^j| < \kappa$, by the regularity of $\kappa$. Finally,
$(\alpha)T^j < \kappa$ ($j \geq 1$) implies that $|(\alpha)\langle
T \rangle| < \kappa$, since $\kappa$ is uncountable.

If there exists $\aleph_0 \leq \lambda < \kappa$ such that
$|(\alpha)T| \leq \lambda$ for all $\alpha \in \Omega$, then each
$\lambda_\beta$ can be taken to be $\lambda$. Let us assume
inductively that for some $j > 1$, $\Sigma = (\alpha)T^{j-1}$ has
cardinality $\leq \lambda^{j-1}$. Then, the above argument shows
that
\begin{xlist}\item
$\displaystyle |(\alpha)T^j| \leq \sum_{\sigma \in \Sigma} \lambda
\leq \sum_{\lambda^{j-1}} \lambda = \lambda^j$.
\end{xlist}
Since $\aleph_0 \leq \lambda$, we have $|\lambda^j| = \lambda$ for
each $j \geq 1$. Therefore $|(\alpha)\langle T \rangle| \leq
\sum_\omega\lambda = \lambda$.
\end{proof}

We conclude the section by noting that $S = \Sym(\Omega)$, the
group of all permutations of $\Omega$, is equivalent to $E =
\E(\Omega)$, with respect to our equivalence relation. This result
is known (cf.\ \cite[Theorem~3.3]{HRH}), but we provide a quick
proof, for the convenience of the reader.

\begin{theorem}[{\rm Howie, Ru\v{s}kuc, and Higgins}]\label{perm-map}
There exist $g_1, g_2 \in E$ such that $E = g_1 S g_2$. In
particular, $E \approx_{\aleph_0} S$.
\end{theorem}

\begin{proof}
Since $\Omega$ is infinite, we can write $\Omega = \bigcup_{\alpha
\in \Omega} \Sigma_\alpha$, where the union is disjoint, and for
each $\alpha \in \Omega$, $|\Sigma_\alpha| = |\Omega|$. Let $g_1
\in E$ be an injective map such that $|\Omega \setminus
(\Omega)g_1| = |\Omega|$, and let $g_2 \in E$ be the map that
takes each $\Sigma_\alpha$ to $\alpha$.

Now, let $f \in E$ be any element. For each $\alpha \in \Omega$,
let $\Delta_\alpha$ denote the preimage of $\alpha$ under $f$. Let
$h \in E$ be an injective self-map that embeds $\Delta_\alpha$ in
$\Sigma_\alpha$, for each $\alpha \in \Omega$, and such that for
some $\alpha \in \Omega$, $|\Sigma_\alpha \setminus
(\Delta_\alpha)h| = |\Omega|$. Then $f = hg_2$. Also, since $g_1$
and $h$ are both injective and $|\Omega \setminus (\Omega)h| =
|\Omega| = |\Omega \setminus (\Omega)g_1|$, there is a permutation
$\bar{h} \in S$ such that for all $\alpha \in \Omega$,
$(\alpha)g_1 \bar{h} = (\alpha)h$. Hence, we have $f = g_1 \bar{h}
g_2 \in g_1 S g_2$.
\end{proof}

\section{The countable case}\label{count section}

From now on we will restrict our attention to the case where
$\Omega$ is countable.  It will often be convenient to assume that
$\Omega = \omega$, the set of natural numbers. However, we will
continue using the symbol $\Omega$; partly in order to distinguish
the role of the set as the domain of our maps from its role as an
indexing set in some of the proofs, and partly because in Section
7 we will be interested in arbitrary orderings of the set. The
symbols $\prec_{\aleph_0,E}$, $\preq_{\aleph_0,E}$, and
$\approx_{\aleph_0,E}$ will henceforth be written simply as
$\prec$, $\preq$, and $\approx$, respectively.

We will say that a set $A$ of disjoint nonempty subsets of
$\Omega$ is a {\it partition} of $\Omega$ if the union of the
members of $A$ is $\Omega$. If $U \subseteq E \ (= \E(\Omega))$
and $A$ is a partition of $\Omega$, let us define $U_{(A)} = \{f
\in U : (\Sigma) f \subseteq \Sigma \ \mathrm{for} \ \mathrm{all}
\ \Sigma \in A\}$. Also, if $U \subseteq E$ and $\Sigma \subseteq
\Omega$, let $U_{(\Sigma)} = \{f \in U : (\alpha) f = \alpha \
\mathrm{for} \ \mathrm{all} \ \alpha \in \Sigma\}$. (The notation
$U_{(A)}$ can be considered an extension of the notation
$U_{(\Sigma)}$.)

The main aim of this section is to show that given two subgroups
$G_1, G_2 \subseteq S \ (= \Sym(\Omega))$ that are closed in the
function topology, we have $G_1 \approx G_2$ if and only if $G_1$
and $G_2$ are equivalent as subgroups, in the sense of~\cite{B&S}.
First, we need two preliminary results.

\begin{proposition}\label{partition monoids}
Let $A$ be a partition of $\Omega$. Then $S_{(A)} \approx
E_{(A)}$.
\end{proposition}

\begin{proof}
If $A$ has an infinite member, then $E_{(A)} \approx E$, by
Lemma~\ref{g,h}, and $S_{(A)} \approx S$, by a similar argument.
Hence, the result follows from Theorem~\ref{perm-map}. Let us,
therefore, assume that all members of $A$ are finite, and let us
write $A = \{A_i : i \in \omega\}$. Further, let $n_i = |A_i|$,
for each $i \in \omega$, and write $A_i = \{a(i,0), a(i,1), \dots,
a(i,n_i-1)\}$. Let $B = \{B_i : i \in \omega\}$ be a partition of
$\Omega$ such that for each $i \in \omega$, $|B_i| = n_i^2$. By
Theorems 13, 15, and 16 of~\cite{B&S}, $S_{(A)} \approx S_{(B)}$.
We will show that $E_{(A)} \preq S_{(B)}$.

For each $i \in \omega$, write $B_i = \{b(i,j,k) : j,k \in \{0, 1,
\dots, n_i-1\}\}$. Let $g_1 \in E$ be the endomorphism that maps
each $A_i$ into $B_i$ via $a(i,j) \mapsto b(i,j,0)$, and let $g_2
\in E$ be the endomorphism that maps each $B_i$ onto $A_i$ via
$b(i,j,k) \mapsto a(i,k)$.

Consider any element $h \in E_{(A)}$, and for each $a(i,j) \in
\Omega$ write $(a(i,j))h = a(i,c_{ij})$, for some $c_{ij} \in \{0,
1, \dots, n_i-1\}$. Let $\bar{h} \in S_{(B)}$ be any permutation
such that for each $i \in \omega$ and $j \in \{0, 1, \dots,
n_i-1\}$, $(b(i,j,0))\bar{h} = b(i,j,c_{ij})$ (e.g., we can define
$\bar{h}$ by $b(i,j,k) \mapsto b(i,j,k+c_{ij} (\mathrm{mod} \,
n_i-1))$). Then, given any $a(i,j) \in \Omega$, we have
$(a(i,j))g_1\bar{h}g_2 = (b(i,j,0))\bar{h}g_2 = (b(i,j,c_{ij}))g_2
= a(i,c_{ij}) = (a(i,j))h$, and hence $E_{(A)} \subseteq g_1
S_{(B)} g_2$.
\end{proof}

In the above proof, we called on results from~\cite{B&S} to deduce
that $S_{(A)} \approx S_{(B)}$. This was done primarily in the
interests of space, since while it is not very difficult to show
that $S_{(A)} \approx S_{(B)}$ directly, several different cases
would need to be considered.

\begin{lemma}\label{partitions & orders}
Let $A$ be a partition of $\Omega$ into finite sets such that
there is no common finite upper bound on the cardinalities of the
members of $A$. Then $E_{(A)} \approx E_\leq$. $($See
Definition~\ref{leq} for the notation $E_\leq$.$)$
\end{lemma}

\begin{proof}
To prove that $E_\leq \preq E_{(A)}$, we will construct $g, h \in
E$ such that $E_\leq \subseteq g E_{(A)}h$. By our hypotheses on
$A$, we can find $\{B_i \in A : i \in \omega\}$, consisting of
disjoint sets, such that each $|B_i| \geq i$. Let us write each
$B_i$ as $\{b(i,0), b(i,1), \dots, b(i,m_i - 1)\}$, where $m_i =
|B_i|$. Define $g \in E$ by $(i)g = b(i,0)$ for all $i \in
\Omega$, and define $h \in E$ by $(b(i,j))h = j$ for all $i \in
\omega$ and $j < m_i$ ($h$ can be defined arbitrarily on elements
not of the form $b(i,j)$). Then $E_\leq \subseteq gE_{(A)}h$.

To prove that $E_{(A)} \preq E_\leq$, write $A = \{A_i : i \in
\omega\}$, and let $g \in E$ be any injective map such that for
all $i \in \omega$ and $a \in A_i$, $(a)g$ is greater than all the
elements of $A_i \ (\subseteq \omega)$. Now, let $f \in E_{(A)}$
be any element. Then for all $a \in \Omega$, $(a)f < (a)g$, since
$a \in A_i$ for some $i \in \omega$. Thus, we can find a map $h
\in E_\leq$ such that for all $a \in \Omega$, $((a)g)h = (a)f$,
since $g$ is injective. Therefore, $E_{(A)} \subseteq gE_\leq$.
\end{proof}

\begin{theorem}\label{monoids&groups}
Let $G_1$ and $G_2$ be subgroups of $S = \Sym(\Omega) \subseteq E$
that are closed in the function topology on $S$. Let us write $G_1
\approx_S G_2$ if $G_1$ and $G_2$ are equivalent as groups
$($i.e., if the group generated by $G_1 \cup U$ is equal to the
group generated by $G_2 \cup U$, for some finite set $U \subseteq
S$$)$. Then $G_1 \approx_S G_2$ if and only if $G_1 \approx G_2$.
\end{theorem}

\begin{proof}
To show the forward implication, suppose that $U \subseteq S$ is a
finite subset such that the group generated by $G_1 \cup U$ is
equal to the group generated by $G_2 \cup U$. Letting $U^{-1}$ be
the set consisting of the inverses of the elements of $U$, we see
that $\langle G_1 \cup (U \cup U^{-1}) \rangle = \langle G_2 \cup
(U \cup U^{-1}) \rangle$. Hence $G_1 \approx_S G_2$ implies $G_1
\approx G_2$.

For the converse, let $A$ be a partition of $\Omega$ into finite
sets such that there is no common finite upper bound on the
cardinalities of the members of $A$, and let $B$ be a partition of
$\Omega$ into 2-element sets. By the main results of~\cite{B&S},
every closed subgroup of $S$ is $\approx_S$ to exactly one of $S$,
$S_{(A)}$, $S_{(B)}$, or $\{1\}$. We finish the proof by showing
that these four groups are $\not\approx$ to each other.

By Proposition~\ref{partition monoids}, $S_{(A)} \approx E_{(A)}$,
and by the previous lemma, the latter is $\approx E_\leq$. By
Theorem~\ref{perm-map}, $S \approx E$. Part (i) of
Theorem~\ref{fear} then implies that $S_{(A)} \prec S$, and part
(ii) of that theorem implies that $S_{(B)} \prec S_{(A)}$. Also,
$\{1\} \prec S_{(B)}$, since $S_{(B)}$ is uncountable. Thus, we
have $\{1\} \prec S_{(B)} \prec S_{(A)} \prec S$.
\end{proof}

\section{An example}

The goal of this section is to show that the partial ordering
$\preq$ of submonoids of $E = \E(\Omega)$ is not a total ordering,
i.e., that there are submonoids $M, M' \subseteq E$ such that $M
\not\preq M'$ and $M' \not\preq M$. In case the reader wishes to
skip this section, we note that nothing in subsequent sections
will depend on the present discussion.

As in Section 1, upon identifying $\Omega$ with $\omega$, let $M_2
\subseteq E$ be the submonoid generated by all maps whose images
are contained in $\{0, 1, 2 \}$. Let $\Sigma_1, \Sigma_2 \subseteq
\Omega$ be disjoint infinite subsets, such that $\{0,1\} \subseteq
\Sigma_1$, $\{2,3\} \subseteq \Sigma_2$, and $\Sigma_1 \cup
\Sigma_2 = \Omega$. Let $M'_3 \subseteq E$ be the submonoid
generated by all maps that take $\Sigma_1$ to $\{0,1\}$ and
$\Sigma_2$ to $\{2,3\}$. We will show that $M_2 \not\preq M'_3$
and $M'_3 \not\preq M_2$.

Using the same argument as in Section 1, it is easy to see that
$M'_3 \not\preq M_2$. (If $g_1 g_2 \dots g_k$ is any word in
elements of $E$, where at least one of the $g_i \in M_2$, then the
image of $\Omega$ under $g_1 g_2 \dots g_k$ has cardinality at
most $3$. Hence, if $M'_3 \subseteq \langle M_2 \cup U \rangle$
for some finite subset $U \subseteq E$, then all the elements of
$M'_3$ whose images have cardinality $4$ must be in $\langle U
\rangle$. This is impossible, since there are uncountably many
such elements.)

Next, let $U \subseteq E$ be a finite set. We will show that $M_2
\not\subseteq \langle M'_3 \cup U \rangle$. We begin by
characterizing the elements of $\langle M'_3 \cup U \rangle$. Let
$H \subseteq E$ be the (countable) set of all maps that fix
$\Omega \setminus \{0,1,2,3\}$ elementwise. Now, consider any word
$f = g_1 g_2 \dots g_k$ in elements of $E$, such that $g_1 \in
M'_3$. Since the image of $g_1$ is contained in $\{0,1,2,3\}$, $f$
can be written as $g_1h$, where $h$ is some element of $H$. Hence,
any element $f \in \langle M'_3 \cup U \rangle$ can be written as
$f = g_1g_2h$, where $g_1 \in \langle U \rangle$, $g_2 \in M'_3$,
and $h \in H$. (Here we are using that fact that $1 \in \langle U
\rangle \cap M'_3 \cap H$.)

We note that each $g \in \langle U \rangle$ either takes
infinitely many elements of $\Omega$ to $\Sigma_1$ or takes
infinitely many elements to $\Sigma_2$, since $\Sigma_1 \cup
\Sigma_2 = \Omega$. For each such $g$, let $\Gamma_g \subseteq
\Omega$ denote either the set of those elements that $g$ takes to
$\Sigma_1$ or the set of those elements that $g$ takes to
$\Sigma_2$ - whichever is infinite. Set $F = \{\Gamma_g : g \in
\langle U \rangle\}$. Since $\langle U \rangle$ is countable, so
is $F$, and hence we can write it as $F = \{\Delta_i : i \in
\omega\}$.

Next, let us construct for each $i \in \omega$ a triplet of
distinct elements $a_i, b_i, c_i \in \Delta_i$, such that the sets
$\{a_i, b_i, c_i\}$ are disjoint. We take $a_0, b_0, c_0 \in
\Delta_0$ to be any three distinct elements (which must exist,
since $\Delta_0$ is infinite). Let $0 \leq j$ be an integer, and
assume that the elements $a_i, b_i, c_i \in \Delta_i$ have been
picked for all $i \leq j$. Let $a_{j+1}, b_{j+1}, c_{j+1} \in
\Delta_{j+1} \setminus \bigcup_{i \leq j} \{a_i, b_i, c_i\}$ be
any three distinct elements. (Again, this is possible, by the fact
that $\Delta_{j+1}$ is infinite.)

Now, let $f \in M_2$ be an element that takes each set $\{a_i,
b_i, c_i\}$ bijectively to $\{0, 1, 2\}$, such that $f \notin
\langle U \rangle$. A self-map with these properties exists, since
there are uncountably many maps that take each $\{a_i, b_i, c_i\}$
bijectively to $\{0, 1, 2\}$, and $\langle U \rangle$ is
countable. We finish the proof by showing that $f \notin \langle
M'_3 \cup U \rangle$. Suppose, on the contrary, that $f \in
\langle M'_3 \cup U \rangle$. Then $f = g_1g_2h$, where $g_1 \in
\langle U \rangle$, $g_2 \in M'_3$, and $h \in H$, by the above
characterization. Since $f \notin \langle U \rangle$, we may
assume that $g_2 \neq 1$. Let $\Delta_i \in F$ be the set
corresponding to $g_1$ (i.e., $\Gamma_{g_1}$). Then, by the above
construction, we can find three distinct elements $a_i, b_i, c_i
\in \Delta_i$ such that $f$ takes $\{a_i, b_i, c_i\}$ bijectively
to $\{0, 1, 2\}$. On the other hand, by choice of $\Delta_i$,
$g_1$ either takes $\{a_i, b_i, c_i\}$ to $\Sigma_1$ or takes
$\{a_i, b_i, c_i\}$ to $\Sigma_2$. In either case, $|(\{a_i, b_i,
c_i\})g_1g_2h| \leq 2$, since $g_2$ takes each of $\Sigma_1$ and
$\Sigma_2$ to a 2-element set. Hence $f \neq g_1g_2h$; a
contradiction. We therefore conclude that $f \notin \langle M'_3
\cup U \rangle$.

In summary, we have

\begin{proposition}
The partial ordering $\preq$ of submonoids of $E$ is not a total
ordering.
\end{proposition}

\section{Four lemmas}

The results of this section (except for the first) are close
analogs of results in~\cite{B&S}. We will use them in later
sections to classify various submonoids of $E = \E(\Omega)$ into
equivalence classes; Lemmas~\ref{tree}\,-\ref{tree3} will be our
main tools for showing that submonoids are $\approx$ to each
other. The proofs of these three lemmas are, for the most part,
simpler than those of their group-theoretic analogs (namely
\cite[Lemma~10]{B&S}, \cite[Lemma~12]{B&S}, and
\cite[Lemma~14]{B&S}, respectively).

\begin{lemma}\label{partition embed}
Let $M \subseteq E$ be a submonoid. Then $M \preq E_{(A)}$, where
$A = \{A_\alpha : \alpha \in \Omega\}$ is any partition of
$\Omega$ such that for each $\alpha \in \Omega$, $|A_\alpha| =
|(\alpha)M|$.
\end{lemma}

\begin{proof}
Let $g \in E$ be a map such that for all $\alpha \in \Omega$,
$(\alpha)g \in A_\alpha$, and let $h \in E$ be a map such that for
all $\alpha \in \Omega$, $h$ maps $(A_\alpha)$ onto $(\alpha)M$.
Then $M \subseteq gE_{(A)}h$.
\end{proof}

\begin{lemma}\label{tree}
Let $M$ be a submonoid of $E$, and suppose there exist a sequence
$(\alpha_i)_{i\in\omega} \in \Omega^\omega$ of distinct elements
and a sequence of nonempty subsets $D_i \subseteq \Omega^i$
$(i\in\omega),$ such that
\begin{enumerate}
\item[\rm{(i)}] For each $i\in\omega$ and
$(\beta_0,\ldots,\beta_{i-1})\in D_i,$ there exist infinitely many
elements $\beta\in\Omega$ such that
$(\beta_0,\ldots,\beta_{i-1},\beta)\in D_{i+1};$ and
\item[\rm{(ii)}] If $(\beta_i)_{i\in\omega}\in\Omega^\omega$ has the property
that $(\beta_0,\ldots,\beta_{i-1})\in D_i$ for each $i\in\omega,$
then there exists $g\in M$ such that for all $i \in \omega$,
$\beta_i=(\alpha_i) g$, and the elements $\beta_i$ are all
distinct.
\end{enumerate}
Then there exist $f, h \in E$ such that $E = fMh$. In particular,
$M \approx E$.
\end{lemma}

\begin{proof}
For each $i \in \omega$ and $(\beta_0,\ldots,\beta_{i-1})\in D_i,$
let
\begin{xlist}\item
$\Gamma(\beta_0,\ldots,\beta_{i-1}) = \{\beta \in \Omega :
(\beta_0,\ldots,\beta_{i-1}, \beta) \in D_{i+1}\}$.
\end{xlist}
By (i), each $\Gamma(\beta_0,\ldots,\beta_{i-1})$ is an infinite
subset of $\Omega$. Since for each $i \in \omega,$ $\Omega^i$ is
countable, $\bigcup_{i \in \omega} \Omega^i$ is countable as well,
and therefore, so is $\bigcup_{i \in \omega} D_i \subseteq
\bigcup_{i \in \omega} \Omega^i$. Thus, there are only countably
many sets of the form $\Gamma(\beta_0,\ldots,\beta_{i-1})$. By a
standard inductive construction, we can find a collection
$\{\Lambda(\beta_0,\ldots,\beta_{i-1}) : i \in \omega, \
(\beta_0,\ldots,\beta_{i-1})\in D_i\}$ of disjoint infinite sets,
such that each $\Lambda(\beta_0,\ldots,\beta_{i-1}) \subseteq
\Gamma(\beta_0,\ldots,\beta_{i-1})$.

Next, we define a map $f \in E$ by
\begin{xlist}\item
$(i)f = \alpha_i$ for all $i \in \omega \ (=\Omega)$.
\end{xlist}
Also, let $h \in E$ be any map that takes each
$\Lambda(\beta_0,\ldots,\beta_{i-1})$ surjectively to $\Omega$. We
will show that $E = fMh$.

Let $g \in E$ be any element. We define recursively a sequence
$(\beta_i)_{i\in\omega} \in \Omega^\omega$ as follows: for each $i
\in \omega$ let $\beta_i \in \Lambda(\beta_0,\ldots,\beta_{i-1})$
be such that $(\beta_i)h = (i)g$. Since each
$\Lambda(\beta_0,\ldots,\beta_{i-1}) \subseteq
\Gamma(\beta_0,\ldots,\beta_{i-1})$, our sequence
$(\beta_i)_{i\in\omega}$ has the property that
$(\beta_0,\ldots,\beta_{i-1})\in D_i$ for each $i\in\omega.$ Thus,
by (ii), there exists $\bar{g} \in M$ such that for all $i \in
\omega$, $\beta_i=(\alpha_i) \bar{g}$. For each $i \in \omega \
(=\Omega)$, we then have that $(i)f\bar{g}h = (\alpha_i)\bar{g}h =
(\beta_i)h = (i)g$, and therefore $g = f\bar{g}h$.
\end{proof}

The next argument uses the same basic idea, but it is more
complicated.

\begin{lemma}\label{tree2}
Let $M$ be a submonoid of $E$, and suppose there exist a sequence
$(\alpha_i)_{i\in\omega}\in\Omega^\omega$ of distinct elements, an
unbounded sequence of positive integers $(N_i)_{i\in\omega},$ and
a sequence of nonempty subsets $D_i \subseteq \Omega^i$
$(i\in\omega),$ such that
\begin{enumerate}
\item[\rm{(i)}] For each $i\in\omega$ and
$(\beta_0,\ldots,\beta_{i-1})\in D_i,$ there exist at least $N_i$
elements $\beta\in\Omega$ such that
$(\beta_0,\ldots,\beta_{i-1},\beta)\in D_{i+1};$ and
\item[\rm{(ii)}] If $(\beta_i)_{i\in\omega}\in\Omega^\omega$ has the property
that $(\beta_0,\ldots,\beta_{i-1})\in D_i$ for each $i\in\omega,$
then there exists $g\in M$ such that for each $i \in \omega$,
$\beta_i=(\alpha_i) g$, and the elements $\beta_i$ are all
distinct.
\end{enumerate}
Then there exist $f, h \in E$ such that $E_\leq \subseteq fMh$. In
particular, $E_\leq \preq M$.
\end{lemma}

\begin{proof}
We will construct recursively integers $i(-1) < i(0) < \ldots <
i(j) < \ldots,$ and for each $j \geq 0$ a subset $C_{i(j)}
\subseteq D_{i(j)}.$

Set $i(-1) = -1$ and $i(0) = 0$, and let $C_{i(0)} = D_0$ be the
singleton consisting of the empty string. Now assume inductively
for some $j \geq 1$ that $i(0), \ldots, i(j{-}1)$ and $C_{i(0)},
\ldots, C_{i(j-1)}$ have been constructed. Let $i(j)$ be an
integer such that $N_{i(j)} > j \cdot |C_{i(j-1)}| +
\sum_{k=0}^{j-1} |C_{i(k)}|$. Let $C_{i(j)} \subseteq D_{i(j)}$ be
a finite subset that for each $(\beta_0, \ldots,
\beta_{i(j-1)-1})\in C_{i(j-1)}$ contains $j$ elements of the form
$(\beta_0, \ldots, \beta_{i(j-1)-1}, \ldots, \beta_{i(j)-2},
\beta)$, such that the elements $\beta$ are distinct from each
other and from all elements that occur as last components of
elements of $C_{i(0)}, \ldots, C_{i(j-1)}$. Our choice of $i(j)$
and condition (i) makes this definition possible. (Actually, it
would have sufficed to pick $i(j)$ so that
$N_{i(j-1)+1}N_{i(j-1)+2} \dots N_{i(j)} > j \cdot |C_{i(j-1)}| +
\sum_{k=0}^{j-1} |C_{i(k)}|$.)

Once the above integers and subsets have been constructed, let us
use the sets $C_{i(j)}$ to construct subsets $F_{i(j)} \subseteq
\Omega^j$. Set $F_{i(0)} = C_{i(0)}$. For each element
$(\beta_{0}, \ldots, \beta_{i(j)-1}) \in C_{i(j)}$ with $j \geq
1$, we define a sequence $(\gamma_{i(0)},\ldots,\gamma_{i(j-1)})$
by setting $\gamma_{i(k)} = \beta_{i(k+1)-1}$ ($0 \leq k \leq
j-1$); i.e., we drop the $\beta_k$ that do not occur as last
components of elements of $C_{i(0)}, \ldots, C_{i(j-1)}$. For each
$j \geq 1$, we then let $F_{i(j)}$ consist of the tuples
$(\gamma_{i(0)}, \ldots, \gamma_{i(j-1)})$. Also, for each element
$(\gamma_{i(0)}, \ldots, \gamma_{i(j-1)}) \in F_{i(j)}$ let
\begin{xlist}\item
$\Gamma(\gamma_{i(0)}, \ldots, \gamma_{i(j-1)}) = \{\beta \in
\Omega : (\gamma_{i(0)}, \ldots, \gamma_{i(j-1)}, \beta) \in
F_{i(j+1)}\}$.
\end{xlist}
By construction, each $|\Gamma(\gamma_{i(0)}, \ldots,
\gamma_{i(j-1)})| \geq j+1$; for simplicity, we will assume that
this is an equality, after discarding some elements if necessary.
Let $h \in E$ be a map such that
\begin{xlist}\item
$h$ takes each $\Gamma(\gamma_{i(0)}, \ldots, \gamma_{i(j-1)})$
onto $\{0, \ldots, j\}$.
\end{xlist}
Such a map exists, since the sets $\Gamma(\gamma_{i(0)}, \ldots,
\gamma_{i(j-1)})$ are all disjoint. Also, let $f \in E$ be defined
by
\begin{xlist}\item
$(j)f = \alpha_{i(j)}$ for all $j \in \omega \ (=\Omega)$.
\end{xlist}
We finish the proof by showing that $E_\leq \subseteq fMh$.

Let $g \in E_\leq$ be any element. We first construct recursively
a sequence $(\gamma_{i(j)})_{j\in\omega}$ such that for each $j
\in \omega$, $(\gamma_{i(0)}, \ldots, \gamma_{i(j-1)}) \in
F_{i(j)}$. Let $\gamma_{i(0)}$ be the unique element of
$\Gamma()$. (We note that $(\gamma_{i(0)})h = 0 = (0)g$, by
definition of $h$.) Assuming that $\gamma_{i(0)}, \ldots,
\gamma_{i(j-1)}$ have been defined, let $\gamma_{i(j)} \in
\Gamma(\gamma_{i(0)}, \ldots, \gamma_{i(j-1)})$ be such that
$(\gamma_{i(j)})h = (j)g$. (Such an element exists, by our
definition of $h$ and the fact that for all $k \in \omega$, $(k)g
\leq k$.) Since the sequence $(\gamma_{i(j)})_{j \in \omega} \in
\Omega^\omega$ has the property that
$(\gamma_{i(0)},\ldots,\gamma_{i(j-1)}) \in F_{i(j)}$ for each $j
\in \omega,$ there exists $\bar{g} \in M$ such that for all $i \in
\omega$, $\gamma_{i(j)}=(\alpha_{i(j)}) \bar{g}$. This follows
from (ii), since $(\gamma_{i(j)})_{j\in\omega}$ is a subsequence
of some $(\beta_i)_{i\in\omega}$ as in (ii). Hence, for each $j
\in \omega \ (=\Omega)$, we have that $(j)f\bar{g}h =
(\alpha_{i(j)})\bar{g}h = (\gamma_{i(j)})h = (j)g$, and therefore
$g = f\bar{g}h$.
\end{proof}

\begin{lemma}\label{tree3}
Let $M$ be a submonoid of $E$. Suppose there exist three sequences
$(\alpha_i)_{i\in\omega},$ $(\beta_i)_{i\in\omega},$
$(\gamma_i)_{i\in\omega} \in \Omega^\omega$ of distinct elements,
such that $(\beta_i)_{i\in\omega}$ and $(\gamma_i)_{i\in\omega}$
are disjoint, and for every element $(\delta_i)_{i\in\omega} \in
\prod_{i\in\omega} \{\beta_i, \gamma_i\} \subseteq \Omega^\omega$,
there exists $g\in M$ such that for all $i \in \omega$, $\delta_i
= (\alpha_i) g$. Then there exist $f, h \in E$ such that $E_{(A)}
\subseteq fMh$, where $A$ is a partition of $\Omega$ into
$2$-element sets. In particular $E_{(A)} \preq M$.
\end{lemma}

\begin{proof}
Write $A = \{A_i : i \in \omega\}$, where for each $i \in \omega,$
$A_i = \{a_i, b_i\}$. Let $f \in E$ be the map defined by $(a_i)f
= \alpha_{2i}$ and $(b_i)f = \alpha_{2i+1}$, and let $h \in E$ be
a map that for each $i \in \omega$ takes $\{\beta_{2i},
\beta_{2i+1}\}$ to $a_i$ and $\{\gamma_{2i}, \gamma_{2i+1}\}$ to
$b_i$. Then $E_{(A)} \subseteq fMh$.
\end{proof}

\section{Submonoids arising from preorders}\label{preordsect}

\begin{definition}
Given a preorder $\rho$ on $\Omega$, let $E(\rho) \subseteq E \ (=
\E(\Omega))$ denote the subset consisting of all maps $f$ such
that for all $\alpha \in \Omega$ one has $(\alpha, (\alpha)f) \in
\rho$.
\end{definition}

Clearly, subsets of the form $E(\rho)$ are submonoids. The
submonoids $E_{(A)}$ (where $A$ is a partition of $\Omega$) are of
this form, as is $E_\leq$. The goal of this section is to classify
such submonoids into equivalence classes. To facilitate the
discussion, let us divide them into five types.

\begin{definition}\label{types}
Let $\rho$ be a preorder on $\Omega$. For each $\alpha \in \Omega$
set $\Delta_\rho (\alpha) = \{\beta \in \Omega : (\alpha, \beta)
\in \rho\}$ $($the ``principal up-set generated by $\alpha$"$)$.
We will say that
\begin{enumerate}
\item[] The preorder $\rho$ is of {\em type 1} if there is an
infinite subset $\Gamma \subseteq \Omega$ such that for all
$\alpha \in \Gamma$, $\Delta_\rho (\alpha)$ is infinite.

\item[] The preorder $\rho$ is of {\em type 2} if the cardinalities of
the sets $\Delta_\rho (\alpha)$ $($$\alpha \in \Omega$$)$ have no
common finite upper bound, but $\Delta_\rho (\alpha)$ is infinite
for only finitely many $\alpha$.

\item[] The preorder $\rho$ is of {\em type 3} if there is a number $n
\in \omega$ such that $|\Delta_\rho (\alpha)| \leq n$ for all but
finitely many $\alpha \in \Omega$, and there are infinitely many
$\alpha \in \Omega$ such that $|\Delta_\rho (\alpha)| > 1$.

\item[] The preorder $\rho$ is of {\em type 4} if $|\Delta_\rho
(\alpha)| = 1$ for all but finitely many $\alpha \in \Omega$.
\end{enumerate}

Let us further divide preorders of type $3$ into two sub-types. We
will say that
\begin{enumerate}
\item[] The preorder $\rho$ is of {\em type 3a} if it is of type
$3$ and, in addition, there are infinite families
$\{\alpha_i\}_{i\in\omega}, \{\beta_i\}_{i\in\omega} \subseteq
\Omega$ consisting of distinct elements, such that for each $i \in
\omega$, $\alpha_i \neq \beta_i$ and $(\alpha_i, \beta_i) \in
\rho$.

\item[] The preorder $\rho$ is of {\em type 3b} if it is of type $3$ and, in
addition, there is a finite set $\Gamma \subseteq \Omega$ such
that for all but finitely many $\alpha \in \Omega$, $\Delta_\rho
(\alpha) \subseteq \Gamma \cup \{\alpha\}$.
\end{enumerate}
\end{definition}

It is clear that every preorder on $\Omega$ falls into exactly one
of the above five types. Further, if $\rho$ is a preorder of type
2 or 3, and $\Sigma = \{\alpha \in \Omega : |\Delta_\rho (\alpha)|
= \aleph_0\}$, then $E(\rho) \approx E(\rho)_{(\Sigma)}$. (The
notation $E(\rho)_{(\Sigma)}$ is defined at the beginning of
Section~\ref{count section}.) For, if $\alpha \in \Omega \setminus
\Sigma$, then $\Delta_\rho (\alpha) \cap \Sigma = \emptyset$,
since $\Delta_\rho (\alpha)$ is finite. This shows that $(\Omega
\setminus \Sigma)E(\rho) \cap \Sigma = \emptyset$. Hence, if $f
\in E(\rho)$ is any element, then $f \in E(\rho)_{(\Sigma)}U$,
where $U = E_{(\Omega \setminus \Sigma)}$. Therefore, $E(\rho)
\subseteq \langle E(\rho)_{(\Sigma)} \cup U \rangle$, and $E(\rho)
\approx E(\rho)_{(\Sigma)}$. (Since $\Sigma$ is finite, $U$ is
countable. Theorem~\ref{2-gen} then allows us to replace $U$ by a
finite set.) We will use this observation a number of times in
this section.

\begin{lemma}\label{preord1}
Let $\rho$ be a preorder on $\Omega$, and let $A$ be a partition
of $\Omega$ into $2$-element sets.
\begin{enumerate}
\item[\rm{(i)}] If $\rho$ is of  type $1$, then $E(\rho) \approx
E$.

\item[\rm{(ii)}] If $\rho$ is of  type $2$, then $E(\rho) \approx
E_\leq$.

\item[\rm{(iii)}] If $\rho$ is of  type $3a$, then $E(\rho)
\approx E_{(A)}$.
\end{enumerate}
\end{lemma}

\begin{proof}
If $\rho$ is of  type 1, respectively type 2, then $E(\rho)$
clearly satisfies the hypotheses of Lemma~\ref{tree}, respectively
Lemma~\ref{tree2}. Hence, $E \approx E(\rho)$, respectively
$E_\leq \preq E(\rho)$. Further, if $\rho$ is of  type 3a, then
upon removing some elements if necessary, we may assume that the
$\{\alpha_i\}_{i\in\omega}$ and $\{\beta_i\}_{i\in\omega}$
provided by the definition are disjoint. $E(\rho)$ then satisfies
the hypotheses of Lemma~\ref{tree3} (with $\alpha_i = \gamma_i$).
Hence $E_{(A)} \preq E(\rho)$.

To finish the proof, we must show that $E(\rho) \preq E_\leq$ if
$\rho$ is of type 2, and that $E(\rho) \preq E_{(A)}$ if $\rho$ is
of type 3a. In either case, we may assume that for all $\alpha \in
\Omega$, $\Delta_\rho (\alpha)$ is finite, by the remarks
following Definition~\ref{types}. Now, $E(\rho) \preq E_{(B)}$,
where $B$ is a partition of $\Omega$ into finite sets, by
Lemma~\ref{partition embed}, and in the case where $\rho$ is of
type 3a, these finite sets can all be taken to have cardinality
$n$, for some $n \in \omega$. If $\rho$ is of type 2, then we have
$E(\rho) \preq E_{(B)} \approx E_\leq$, by Lemma~\ref{partitions &
orders}. If $\rho$ is of type 3a, then $E_{(A)} \approx E_{(B)}$,
by Proposition~\ref{partition monoids} and \cite[Theorem~15]{B&S}.
Hence $E(\rho) \preq E_{(A)}$.
\end{proof}

\begin{definition}
Given $\gamma \in \Omega$, let $\rho_\gamma$ be the preorder on
$\Omega$ defined by $(\alpha, \beta) \in \rho_\gamma
\Leftrightarrow \beta \in \{\alpha, \gamma\}$.
\end{definition}

It is clear that the preorders $\rho_\gamma$ are of type $3b$.
Further, if $\alpha, \beta \in \Omega$ are any two elements, then
$E(\rho_\alpha) \approx E(\rho_\beta)$. More specifically, if $g
\in E$ is any permutation of $\Omega$ that takes $\alpha$ to
$\beta$, then $E(\rho_\alpha) = gE(\rho_\beta)g^{-1}$. (Given
$\gamma \in \Omega$, an element of $gE(\rho_\beta)g^{-1}$ either
fixes $\gamma$ or takes it to $\alpha$. Hence
$gE(\rho_\beta)g^{-1} \subseteq E(\rho_\alpha)$, and similarly
$g^{-1}E(\rho_\alpha)g \subseteq E(\rho_\beta)$. Conjugating the
latter expression by $g$, we obtain $E(\rho_\alpha) \subseteq
gE(\rho_\beta)g^{-1}$.)

Let $E_{3b}$ denote the monoid generated by $\{E(\rho_\alpha) :
\alpha \in \Omega \}$, and let $g \in E$ be a permutation which is
transitive on $\Omega$. Then, by the previous paragraph, $E_{3b}
\subseteq \langle E(\rho_\gamma) \cup \{g\} \rangle$ for any
$\gamma \in E$. Hence, for all $\gamma \in E$, $E(\rho_\gamma)
\approx E_{3b}$. As an aside, we note that $E_{3b}$ is closed
under conjugation by permutations of $\Omega$. (Given any
permutation $g$ and any $\alpha \in \Omega$, $gE(\rho_\beta)g^{-1}
= E(\rho_\alpha)$, where $\beta = (\alpha)g$. Thus $gE_{3b}g^{-1}$
contains all the generators of $E_{3b}$, and therefore $E_{3b}
\subseteq gE_{3b}g^{-1}$. Conjugating by $g^{-1}$, we obtain
$g^{-1}E_{3b}g \subseteq E_{3b}$.)

\begin{lemma}\label{preord2}
Let $\rho$ and $\rho'$ be preorders on $\Omega$ of type $3b$. Then
$E(\rho) \approx E(\rho')$.
\end{lemma}

\begin{proof}
Let $\rho$ be a preorder on $\Omega$ of type $3b$. Then there is a
$\beta \in \Omega$ and an infinite set $\Sigma \subseteq \Omega$
such that all $\alpha \in \Sigma$, we have $(\alpha, \beta) \in
\rho$. Let $f_1 \in E$ be a map that takes $\Omega$ bijectively to
$\Sigma$, while fixing $\beta$ (which, we may assume, is an
element of $\Sigma$), and let $f_2 \in E$ be a right inverse of
$f_1$. Then $E(\rho_\beta) \subseteq f_1 E(\rho) f_2$, and hence
$E_{3b} \preq E(\rho)$. We conclude the proof by showing that
$E(\rho) \preq E_{3b}$.

By the remarks following Definition~\ref{types}, we may assume
that for all $\alpha \in \Omega$, $\Delta_\rho (\alpha)$ is
finite. Let $\Gamma \subseteq \Omega$ be a finite set such that
for all $\alpha \in \Omega$, $\Delta_\rho (\alpha) \subseteq
\Gamma \cup \{\alpha\}$, and write $\Gamma = \{\alpha_0, \alpha_1,
\dots, \alpha_n\}$. We will first show that $E(\rho)_{(\Gamma)}
\subseteq E_{3b}$.

Let $h \in E(\rho)_{(\Gamma)}$ be any element. Then we can write
$\Omega$ as a disjoint union $\Lambda_0 \cup \Lambda_1 \cup \dots
\cup \Lambda_n \cup \Lambda$, where for all $\beta \in \Lambda_i$,
$(\beta)h = \alpha_i$, and $h$ acts as the identity on $\Lambda$.
For each $i \in \{0, 1, \dots, n\}$, let $g_i \in
E(\rho_{\alpha_i})$ be the map that takes $\Lambda_i \setminus
\{\alpha_0, \alpha_1, \dots, \alpha_n\}$ to $\alpha_i$ and fixes
all other elements. Then $h = g_0g_1 \dots g_n \in E_{3b}$, and
hence $E(\rho)_{(\Gamma)} \subseteq E_{3b}$.

For each $g \in \E(\Gamma)$ let $f_g \in E$ be such that $f_g$
acts as $g$ on $\Gamma$ and as the identity elsewhere, and set $V
= \{f_g : g \in \E(\Gamma)\}$. Now, let $h \in E(\rho)$ be any
element, and define $\bar{h} \in E(\rho)_{(\Gamma)}$ by $(\alpha)
\bar{h} = (\alpha)h$ for all $\alpha \notin \Gamma$. Noting that,
by definition of $\Gamma$, $(\Gamma)h \subseteq \Gamma$, there is
an element $f_g \in V$ that agrees with $h$ on $\Gamma$. Then $h =
f_g \bar{h} \in VE(\rho)_{(\Gamma)} \subseteq V E_{3b}$, and hence
$E(\rho) \preq E_{3b}$.
\end{proof}

We will say that a map $f \in E$ is {\em finitely-many-to-one}
(or, more succinctly, {\em fm-to-one}) if the preimage of each
element of $\Omega$ under $f$ is finite.

\begin{lemma}\label{preord3}
Let $\rho$ and $\rho'$ be preorders on $\Omega$ of types $3a$ and
$3b$, respectively. Then $E(\rho') \prec E(\rho)$.
\end{lemma}

\begin{proof}
Let $A$ be a partition of $\Omega$ into $2$-element sets, and let
us fix an element $\gamma \in \Omega$. By the previous two lemmas,
it suffices to show that $E(\rho_\gamma) \prec E_{(A)}$. We begin
by proving that $E_{(A)} \not\preq E(\rho_\gamma)$.

For each finite subset $\Sigma \subseteq \Omega$ let $f_\Sigma \in
E$ be the map defined by
\begin{xlist}\item
$(\alpha) f_\Sigma = \left\{ \begin{array}{cl}
\gamma & \textrm{if } \alpha \in \Sigma\\
\alpha & \textrm{if } \alpha \notin \Sigma.\\
\end{array}\right.$
\end{xlist}
Now, let $U \subseteq E$ be a finite subset, and consider a monoid
word $g = g_0 g_1 \dots g_{i-1} g_i g_{i+1} \dots g_n$ in elements
of $E(\rho_\gamma) \cup U$. Suppose that $g$ is fm-to-one and that
the element $g_i$ is in $E(\rho_\gamma)$. Let $\Gamma \subseteq
\Omega$ be the preimage of $\gamma$ under $g_i$. Then $\Sigma :=
((\Omega)g_0 \dots g_{i-1}) \cap \Gamma$ must be finite, and so we
have $g = g_0 \dots g_{i-1} f_\Sigma g_{i+1} \dots g_n$. In a
similar fashion, assuming that $g$ is fm-to-one, we can replace
every element of $E(\rho_\gamma)$ occurring in the word $g$ by an
element of the form $f_\Sigma$, for some finite $\Sigma \subseteq
\Omega$. Considering that all elements of $E_{(A)}$ are fm-to-one,
we conclude that if $h \in E_{(A)} \cap \langle E(\rho_\gamma)
\cup U \rangle$, then $h \in \langle \{f_\Sigma : \Sigma \subseteq
\Omega \ \mathrm{finite} \}\cup U \rangle$. But, the latter set is
countable and hence cannot contain all of $E_{(A)}$. Therefore
$E_{(A)} \not\preq E(\rho_\gamma)$.

It remains to show that $E(\rho_\gamma) \preq E_{(A)}$. Write $A =
\{A_i : i \in \Omega \}$, and for each $i \in \Omega$ set $A_i =
\{\alpha_{i1}, \alpha_{i2}\}$. Let $g_1 \in E$ be the map defined
by $(i)g_1 = \alpha_{i1}$ ($i \in \Omega$), and let $g_2 \in E$ be
defined by $(\alpha_{i1})g_2 = \alpha_{i1}$ and $(\alpha_{i2})g_2
= \gamma$. Then $E(\rho_\gamma) \subseteq g_1 E_{(A)} g_2$.
\end{proof}

\begin{theorem}\label{partition classification}
Let $\rho$ and $\rho'$ be preorders on $\Omega$. Then $E(\rho)
\approx E(\rho')$ if and only if $\rho$ and $\rho'$ are of the
same type.
\end{theorem}

\begin{proof}
The above lemmas, in conjunction with Theorem~\ref{fear}, give the
desired conclusion if $\rho$ and $\rho'$ are each of type 1, 2,
3a, or 3b. The result then follows from the fact that $\rho$ is of
type 4 if and only if $E(\rho)$ is countable.
\end{proof}

\section{The function topology}\label{topology section}

From now on we will be concerned with submonoids that are closed
in the function topology on $E$, so let us recall some facts about
this topology.

Regarding the infinite set $\Omega$ as a discrete topological
space, the monoid $E = \E(\Omega)$, viewed as the set of all
functions from $\Omega$ to $\Omega$, becomes a topological space
under the function topology. A subbasis of open sets in this
topology is given by the sets $\{f \in E : (\alpha)f = \beta\}$
$(\alpha,\beta \in \Omega).$ The closure of a set $U \subseteq E$
consists of all maps $f$ such that, for every finite subset
$\Gamma \subseteq \Omega,$ there exists an element of $U$ agreeing
with $f$ at all members of $\Gamma.$ It is easy to see that
composition of maps is continuous in this topology. Given a subset
$U \subseteq E,$ we will write $\mathrm{cl}_E (U)$ for the closure
of $U$ in $E$.

The following lemma will be useful later on. It is an analog
of~\cite[Lemma~8]{B&S}, with monoids in place of groups and
forward orbits in place of orbits. (Given an element $\alpha \in
\Omega$ and a subset $U \subseteq E$, we will refer to the set
$(\alpha)U$ as the {\em forward orbit} of $\alpha$ under $U$.) The
proof is carried over from \cite{B&S} almost verbatim.

\begin{lemma}\label{function top}
Let $M \subseteq E$ be a submonoid. Then

\begin{enumerate}
\item[\rm{(i)}] $\mathrm{cl}_E (M)$ is also a submonoid of $E$.
\item[\rm{(ii)}] $M$ and $\mathrm{cl}_E (M)$ have the same forward orbits in
$\Omega$.
\item[\rm{(iii)}] If $\Gamma$ is a finite subset of $\Omega,$ then
$\mathrm{cl}_E (M)_{(\Gamma)} = \mathrm{cl}_E (M_{(\Gamma)}).$
\end{enumerate}
\end{lemma}

\begin{proof}
Statement (i) follows from the fact that composition of maps is
continuous.

From the characterization of the closure of a set in our topology,
we see that for $\alpha,\beta \in \Omega,$ the set $\mathrm{cl}_E
(M)$ will contain elements taking $\alpha$ to $\beta$ if and only
if $M$ does, establishing (ii).

Given any subset $\Gamma \subseteq \Omega$, the elements of
$\mathrm{cl}_E (M_{(\Gamma)})$ fix $\Gamma$ elementwise, by (ii).
Hence, $\mathrm{cl}_E (M)_{(\Gamma)} \supseteq \mathrm{cl}_E
(M_{(\Gamma)})$. To show $\mathrm{cl}_E (M)_{(\Gamma)} \subseteq
\mathrm{cl}_E (M_{(\Gamma)})$, assume that $\Gamma \subseteq
\Omega$ is finite, and let $f \in \mathrm{cl}_E (M)_{(\Gamma)}$.
Then every neighborhood of $f$ contains elements of $M$, since $f
\in \mathrm{cl}_E (M)$. But, since $f$ fixes all points of the
finite set $\Gamma$, every sufficiently small neighborhood of $f$
consists of elements which do the same. Hence, every such
neighborhood contains elements of $M_{(\Gamma)}$, and so $f \in
\mathrm{cl}_E (M_{(\Gamma)})$.
\end{proof}

\section{Large stabilizers}

Let us say that a submonoid $M$ of $E = \E(\Omega)$ has {\em large
stabilizers} if for each finite subset $\Sigma \subseteq \Omega$,
$M_{(\Sigma)} \approx M$. For example, all subgroups of $S$ have
large stabilizers, by \cite[Lemma~2]{B&S}, as do submonoids of the
form $E_{(A)}$, where $A$ is a partition of $\Omega$. More
generally, we have

\begin{proposition}\label{large stab ex}
Let $\rho$ be a preorder on $\Omega$. Then $E(\rho)$ has large
stabilizers.
\end{proposition}

\begin{proof}
Let $\Sigma \subseteq \Omega$ be finite. Then $E(\rho)_{(\Sigma)}$
is still a submonoid of the form $E(\rho')$, where $\rho'$ is a
preorder on $\Omega$ of the same type as $\rho$. The desired
conclusion then follows from Theorem~\ref{partition
classification}.
\end{proof}

The following lemma gives another class of submonoids that have
large stabilizers. The proof is similar to the one for
\cite[Lemma~2]{B&S}

\begin{lemma}
Let $G \subseteq \Sym(\Omega) \subseteq E$ be a subgroup. Then
$\mathrm{cl}_E (G)$ has large stabilizers.
\end{lemma}

\begin{proof}
Let $\Sigma \subseteq \Omega$ be finite, and take $f \in
\mathrm{cl}_E (G)$. Then there is a sequence of elements of $G$
that has limit $f$. Eventually, the elements of this sequence must
agree on the members of $\Sigma$, and hence, must lie in some
right coset of $G_{(\Sigma)}$, say the right coset represented by
$g \in G$. Then $f \in \mathrm{cl}_E (G_{(\Sigma)})g =
\mathrm{cl}_E (G)_{(\Sigma)}g$, by Lemma~\ref{function top}.
Letting $R \subseteq G$ be a set of representatives of the right
cosets of $G_{(\Sigma)}$, we have $\mathrm{cl}_E (G) \subseteq
\langle \mathrm{cl}_E (G)_{(\Sigma)} \cup R \rangle$. Now, $R$ is
countable, since $\Omega$ is countable and $\Sigma$ is finite, and
hence $\mathrm{cl}_E (G) \approx \mathrm{cl}_E (G)_{(\Sigma)}$, by
Theorem~\ref{2-gen}.
\end{proof}

Not all submonoids of $E$, however, have large stabilizers. Given
a natural number $n$, the monoid $M_n$ (generated by all maps
whose images are contained in $\{0, 1, \dots, n\}$) mentioned in
Section 1 is an easy example of this. (For any finite $\Sigma
\subseteq \Omega \ (=\omega)$ such that $\Sigma \cap \{0, 1,
\dots, n\} = \emptyset$, $(M_n)_{(\Sigma)} = \{1\} \not\approx
M_n$.)

\begin{definition}
Let us say that a map $f \in E$ is {\em a.e.\ injective} if there
exists some finite set $\Gamma \subseteq \Omega$ such that $f$ is
injective on $\Omega \setminus \Gamma$.
\end{definition}

In much of the sequel we will be concerned with submonoids where
the a.e.\ injective maps form dense subsets (viewing $E$ as
topological space under the function topology). A submonoid $M
\subseteq E$ has this property if and only if for every finite
$\Sigma \subseteq \Omega$ and every $f \in M$, there is an a.e.\
injective map $g \in M$ that agrees with $f$ on $\Sigma$. Clearly,
a.e.\ injective maps form dense subsets in submonoids of $E$ that
consist of injective maps, such as subgroups of $\Sym(\Omega)$ and
their closures in the function topology. It is easy to see that
this is also the case for the submonoids $E(\rho)$, since given
any $g \in E(\rho)$ and any finite $\Sigma \subseteq \Omega$, we
can find an $f \in E(\rho)$ that agrees with $g$ on $\Sigma$ and
acts as the identity elsewhere. From now on we will abuse language
by referring to submonoids in which the a.e.\ injective maps form
dense subsets as submonoids having dense a.e.\ injective maps.

Our goal will be to classify into equivalence classes submonoids
of $E$ that are closed in the function topology, and have large
stabilizers and dense a.e.\ injective maps. By the remarks above,
examples of such submonoids include closures in the function
topology in $E$ of subgroups of $\Sym(\Omega)$ and submonoids of
the form $E(\rho)$. These equivalence classes will be shown to be
represented by $E$, $E_\leq$, $E_{(A)}$, $E(\rho_\gamma)$, and
$\{1\}$, respectively (where $A$ is a partition of $\Omega$ into
2-element sets, and $\gamma \in \Omega$). Most of the work will go
into showing that a given monoid from the above list is $\preq M$,
for a monoid $M$ satisfying an appropriate condition. In the first
three cases (and to some extent in the fourth) our arguments will
follow a similar pattern. We will construct sequences of the form
$g_0, g_0g_1, g_0g_1g_2, \dots$, and use closure in the function
topology to conclude that such sequences have limits in the
appropriate monoid $M$. The subsets consisting of these limits
will then satisfy the hypotheses of one of the
Lemmas~\ref{tree}\,-\ref{tree3}, giving us the desired conclusion.

\begin{proposition}\label{inf orbits}
Let $M \subseteq E$ be a submonoid that is closed in the function
topology and has dense fm-to-one $($and hence a.e.\ injective$)$
maps. If $M_{(\Sigma)}$ has an infinite forward orbit for every
finite $\Sigma \subseteq \Omega$, then $M \approx E$.
\end{proposition}

\begin{proof}
This proof closely follows that of \cite[Theorem~11]{B&S}. We
begin by recursively constructing for each $j \geq 0$ an element
$\alpha_j \in \Omega$ and a finite subset $K_j \subseteq M$,
consisting of fm-to-one maps (``fm-to-one" is defined directly
preceding Lemma~\ref{preord3}), indexed
\begin{xlist}\item
$K_j = \{g(k_0, k_1, \ldots, k_{r-1}) : k_0, k_1, \ldots, k_{r-1},
r \in \omega,\, r + k_0 + \ldots + k_{r-1} = j\}.$
\end{xlist}
Let $\alpha_0$ be any element that has an infinite forward orbit
under $M$. Assume inductively that $\alpha_0, \ldots,
\alpha_{j-1}$ have been defined, and let $\Gamma_j = \{\alpha_0,
\ldots, \alpha_{j-1}\} \cup \{0, \ldots, j-1\}$. We then take
$\alpha_j \in \Omega$ to be any element that has an infinite
forward orbit under $M_{(\Gamma_j)}$.

Now we construct the sets $K_j$. If $j=0$, we have only one
element to choose, $g()$, and we take this to be the identity
element $1 \in M$. Assume inductively that the sets $K_i$ have
been defined for all nonnegative $i < j$. Let us fix arbitrarily
an order in which the elements of $K_j$ are to be constructed.
When it is time to define $g(k_0, k_1, \ldots, k_{r-1})$, let us
write $g' = g(k_0, k_1, \ldots, k_{r-2})$, noting that this is a
member of $K_{j-k_{r-1}-1}$ and hence already defined (and
fm-to-one). We set $g(k_0, k_1, \ldots, k_{r-1}) = hg'$, where $h
\in M_{(\Gamma_{r-1})}$ is chosen so that the image of
$\alpha_{r-1}$ under $hg'$ is distinct from the images of
$\alpha_0, \dots, \alpha_{r-2}$ under the finitely many elements
of $K_0 \cup \ldots \cup K_{j-1}$, and also under the elements of
$K_j$ that have been constructed so far. This is possible, since
$(\alpha_{r-1}) M_{(\Gamma_{r-1})}$ is infinite, by the choice of
$\alpha_{r-1}$, and $(\alpha_{r-1}) M_{(\Gamma_{r-1})}g'$ is
infinite, since $g'$ is fm-to-one. Further, $h$ can be chosen to
be fm-to-one, by our hypothesis on $M$, making $g(k_0, k_1,
\ldots, k_{r-1})$ fm-to-one as well. We note that the images of
$\alpha_0,\ldots,\alpha_{r-2}$ and $0, \ldots, r-2$ under $hg'$
will be the same as their images under $g'$, since elements of
$M_{(\Gamma_{r-1})}$ fix $\{\alpha_0, \ldots, \alpha_{r-2}\} \cup
\{0, \ldots, r-2\} = \Gamma_{r-1}$.

Once the elements of each set $K_j$ are constructed, we have
monoid elements $g(k_0, \ldots, k_{i-1})$ for all $i,\ k_0,
\ldots, k_{i-1} \in \omega.$ We can thus define, for each $i \in
\omega,$
\begin{xlist}\item
$D_i=\{((\alpha_0) g, \ldots, (\alpha_{i-1}) g): g = g(k_0,
\ldots, k_{i-1}) \ \mathrm{for} \ \mathrm{some} \ k_0, \ldots,
k_{i-1} \in \omega\}.$
\end{xlist}
Any two elements of the form $g(k_0, \ldots, k_i)$ with indices
$k_0, \ldots, k_{i-1}$ the same, but different last indices $k_i$,
act differently on $\alpha_i$, so the sets $D_i$ satisfy condition
(i) of Lemma~\ref{tree}. Suppose that $(\beta_i) \in
\Omega^\omega$ has the property that for every $i$ the sequence
$(\beta_0, \ldots, \beta_{i-1})$ is in $D_i$. By construction, the
elements $\beta_i$ are all distinct. Also, we see inductively that
successive strings $(), (\beta_0), \dots,\,(\beta_0, \ldots,
\beta_{i-1}), \ldots$ must arise from unique elements of the forms
$g()$, $g(k_0)$, $\ldots \,$, $g(k_0, \ldots, k_{i-1})$, $\ldots
\,$. By including $\{0, \ldots, r-1\}$ in $\Gamma_{r-1}$ above, we
have ensured that the elements of this sequence agree on larger
and larger subsets of $\omega = \Omega$, which have union
$\Omega$. Thus, this sequence converges to a map $g \in E$, which
necessarily sends $(\alpha_i)$ to $(\beta_i)$. Since $M$ is
closed, $g \in M$; which establishes condition (ii) of
Lemma~\ref{tree}. Hence, that lemma tells us that $M \approx E$.
\end{proof}

In the above proposition, the hypothesis that $M_{(\Sigma)}$ has
an infinite forward orbit for every finite $\Sigma \subseteq
\Omega$ is not necessary for $M \approx E$. For example, let $E_>
\subseteq E$ denote the submonoid generated by maps that are
strictly increasing with respect to the usual ordering of $\omega
= \Omega$ (i.e., maps $f \in E$ such that for all $\alpha \in
\Omega$, $(\alpha)f > \alpha$). This submonoid is closed in the
function topology, and its a.e.\ injective maps form a dense
subset. It is also easy to see that the submonoid $E_>$ satisfies
the hypotheses of Lemma~\ref{tree}, and hence is $\approx E$, but
given any finite $\Sigma \subseteq \Omega$, $(E_>)_{(\Sigma)} =
\{1\}$. However, for submonoids $M \subseteq E$ that have large
stabilizers, the condition that $M_{(\Sigma)}$ has an infinite
forward orbit for every finite $\Sigma \subseteq \Omega$ is
necessary for $M \approx E$.

\begin{corollary}\label{large stab 1}
Let $M \subseteq E$ be a submonoid that is closed in the function
topology, and has large stabilizers and dense fm-to-one maps. Then
$M \approx E$ if and only if $M_{(\Sigma)}$ has an infinite
forward orbit for every finite $\Sigma \subseteq \Omega$.
\end{corollary}

\begin{proof}
If there is a finite set $\Sigma \subseteq \Omega$ such that
$M_{(\Sigma)}$ has no infinite forward orbits, then $M_{(\Sigma)}
\prec E$, by Theorem~\ref{fear}(i). Since $M$ has large
stabilizers, this implies that $M \prec E$. The converse follows
from the previous proposition.
\end{proof}

Let us now turn to submonoids whose stabilizers have finite
forward orbits. We will say that a fm-to-one map $f \in E$ is {\em
boundedly finitely-many-to-one} (abbreviated {\em bfm-to-one}) if
there is a common finite upper bound on the cardinalities of the
preimages of elements of $\Omega$ under $f$.

\begin{proposition}\label{unbounded orbits}
Let $M \subseteq E$ be a submonoid that is closed in the function
topology and has dense bfm-to-one $($and hence a.e.\ injective$)$
maps. If for every finite $\Sigma \subseteq \Omega$ the
cardinalities of the forward orbits of $M_{(\Sigma)}$ have no
common finite bound, then $E_\leq \preq M$.
\end{proposition}

\begin{proof}
Again, this proof closely follows that of \cite[Theorem~13]{B&S}.
Let us fix an unbounded sequence of positive integers
$(N_i)_{i\in\omega}$. We begin by recursively constructing for
each $j \geq 0$ an element $\alpha_j \in \Omega$ and a finite
subset $K_j \subseteq M$, consisting of bfm-to-one maps, indexed
\begin{xlist}\item
$K_j = \{g(k_0, k_1, \ldots, k_{j-1}) : 0 \leq k_i < N_i \ (0 \leq
i < j) \}.$
\end{xlist}
Let the 1-element set $K_0 = \{g()\}$ consist of the identity map
$1 \in M$. Assume inductively that for some $j \geq 0$ the
elements $\alpha_0, \dots, \alpha_{j-1}$ and the sets $K_0, \dots,
K_j$ have been defined, and let $\Gamma_j = \{\alpha_0, \ldots,
\alpha_{j-1}\} \cup \{0, \ldots, j-1\}$. Now, let us choose
$\alpha_j \in \Omega$ so that for each $g = g(k_0, k_1, \ldots,
k_{j-1}) \in K_j$, the forward orbit of $\alpha_j$ under
$M_{(\Gamma_j)}g$ has cardinality at least
\begin{xlist}\item
$\displaystyle (j+1) \cdot (\sum_{i=0}^{j+1}|K_i|) =  (j+1) \cdot
(\sum_{i=0}^j|K_i| + N_0N_1 \dots N_j).$
\end{xlist}
This is possible, since each such $g$ is bfm-to-one, and since
$K_j$ is finite.

Next, let us fix arbitrarily an order in which the elements of
$K_{j+1}$ are to be constructed. When it is time to construct
$g(k_0, k_1, \ldots, k_j)$, let us write $g' = g(k_0, k_1, \ldots,
k_{j-1}) \in K_j$. We define $g(k_0, k_1, \ldots, k_j) = hg'$,
where $h \in M_{(\Gamma_j)}$ is chosen so that the image of
$\alpha_j$ under $hg'$ is distinct from the images of $\alpha_0,
\dots, \alpha_{j-1}$ under the elements of $K_0 \cup \ldots \cup
K_j$ and also under the elements of $K_{j+1}$ that have been
constructed so far. Our choice of $\alpha_j$ makes this possible.
Further, $h$ can be chosen to be bfm-to-one, by our hypothesis on
$M$, making $g(k_0, k_1, \ldots, k_j)$ bfm-to-one as well. We note
that the images of $\alpha_0, \ldots, \alpha_{j-1}$ and $0,
\ldots, j-1$ under $hg'$ will be the same as their images under
$g'$, since elements of $M_{(\Gamma_j)}$ fix $\{\alpha_0, \ldots,
\alpha_{j-1}\} \cup \{0, \ldots, j-1\} = \Gamma_j$.

Once the sets $K_j$ are constructed, for each $i \in \omega,$ let
\begin{xlist}\item
$D_i=\{((\alpha_0) g, (\alpha_1) g, \ldots, (\alpha_{i-1}) g): g
\in K_i\}.$
\end{xlist}
Any two elements of the form $g(k_0, \ldots, k_i)$ with indices
$k_0, \ldots, k_{i-1}$ the same but different last indices $k_i$
act differently on $\alpha_i$, so the sets $D_i$ satisfy condition
(i) of Lemma~\ref{tree2}. Suppose that $(\beta_i) \in
\Omega^\omega$ has the property that for every $i$ the sequence
$(\beta_0, \ldots, \beta_{i-1})$ is in $D_i$. By construction, the
elements $\beta_i$ are all distinct. Also, we see inductively that
successive strings $(), (\beta_0), \dots,\,(\beta_0, \ldots,
\beta_{i-1}), \ldots \,$ must arise from unique elements of the
forms $g()$, $g(k_0)$, $\ldots \,$, $g(k_0, \ldots, k_{i-1})$,
$\ldots \,$. The elements of the above sequence agree on the
successive sets $\{0, \ldots, j-1\}$ (having union $\omega =
\Omega$), and hence the sequence converges to a map $g \in E$,
which must send $(\alpha_i)$ to $(\beta_i)$. Since $M$ is closed,
we have that $g \in M$, establishing condition (ii) of
Lemma~\ref{tree2}. Hence, that lemma tells us that $E_\leq \preq
M$.
\end{proof}

As with Proposition~\ref{inf orbits}, the hypothesis that for
every finite $\Sigma \subseteq \Omega$ the cardinalities of the
forward orbits of $M_{(\Sigma)}$ have no common finite bound is
not necessary for $E_\leq \preq M$. For instance, the submonoid
$E_< \subseteq E$ generated by maps that are strictly decreasing
with respect to the usual ordering of $\omega = \Omega$ (i.e.,
maps $f \in E$ such that for all $\alpha \in \Omega \setminus
\{0\}$, $(\alpha) f < \alpha$, and $(0)f = 0$) is closed in the
function topology, has dense a.e.\ injective maps, and satisfies
the hypotheses of Lemma~\ref{tree2}. Therefore $E_\leq \preq E_<$,
but given any finite $\Sigma \subseteq \Omega$, $(E_<)_{(\Sigma)}
= \{1\}$.

\begin{corollary}\label{large stab 2}
Let $M \subseteq E$ be a submonoid that is closed in the function
topology, and has large stabilizers and dense bfm-to-one maps.
Then $M \approx E_\leq$ if and only if for every finite $\Sigma
\subseteq \Omega$ the cardinalities of the forward orbits of
$M_{(\Sigma)}$ have no common finite bound, and there exists a
finite set $\Sigma \subseteq \Omega$ such that all forward orbits
of $M_{(\Sigma)}$ are finite.
\end{corollary}

\begin{proof}
If $M \approx E_\leq$, then for every finite $\Sigma \subseteq
\Omega$, the cardinalities of the forward orbits of $M_{(\Sigma)}$
have no common finite bound. For, otherwise, there would be some
finite $\Sigma \subseteq \Omega$ for which $M_{(\Sigma)}$ would
satisfy the hypotheses of Theorem~\ref{fear}(ii), implying that $M
\approx M_{(\Sigma)} \prec E_\leq$. Also, there must be a finite
set $\Sigma \subseteq \Omega$ such that all forward orbits of
$M_{(\Sigma)}$ are finite, since otherwise we would have $M
\approx E$, by Proposition~\ref{inf orbits}, contradicting
Theorem~\ref{fear}(i).

Conversely, if for every finite $\Sigma \subseteq \Omega$, the
cardinalities of the forward orbits of $M_{(\Sigma)}$ have no
common finite bound, then $E_\leq \preq M$, by the previous
proposition. Suppose that, in addition, there exists a finite set
$\Sigma \subseteq \Omega$ such that all forward orbits of
$M_{(\Sigma)}$ are finite. Then $M_{(\Sigma)} \preq E_{(A)}$,
where $A$ is a partition of $\Omega$ into finite sets, by
Lemma~\ref{partition embed}. By the large stabilizer hypothesis
and Lemma~\ref{partitions & orders}, we have $M \approx
M_{(\Sigma)} \preq E_{(A)} \approx E_\leq$. Hence $M \approx
E_\leq$.
\end{proof}

The rest of the section is devoted to the more intricate case of
submonoids with stabilizers whose forward orbits have a common
finite bound.

\begin{lemma}
Let $M \subseteq E$ be a submonoid that has large stabilizers, and
let $A$ be a partition of $\Omega$ into $2$-element sets. If there
exists a finite set $\Sigma \subseteq \Omega$ and a positive
integer $n$ such that the cardinalities of all forward orbits of
$M_{(\Sigma)}$ are bounded by $n$, then $M \preq E_{(A)}$.
\end{lemma}

\begin{proof}
Given a finite set $\Sigma \subseteq \Omega$ and a positive
integer $n$ as above, $M_{(\Sigma)} \preq E_{(B)}$, where $B$ is a
partition of $\Omega$ into sets of cardinality $n$, by
Lemma~\ref{partition embed}. But, $E_{(A)} \approx E_{(B)}$, by
Proposition~\ref{partition monoids} and \cite[Theorem~15]{B&S}.
Hence $M \approx M_{(\Sigma)} \preq E_{(A)}$.
\end{proof}

The following proof is based on that of \cite[Theorem~15]{B&S},
but it is more complicated.

\begin{lemma}
Let $M \subseteq E$ be a submonoid that is closed in the function
topology and has dense a.e.\ injective maps. If for all finite
$\Gamma, \Delta \subseteq \Omega$ there exists $\alpha \in \Omega$
such that $(\alpha)M_{(\Delta)} \not\subseteq \Gamma \cup
\{\alpha\}$, then $E_{(A)} \preq M$ $($where $A$ is a partition of
$\Omega$ into $2$-element sets$)$.
\end{lemma}

\begin{proof}
We may assume that there exist finite sets $\Gamma, \Delta
\subseteq \Omega$ and a positive integer $n$ such that for all
$\alpha \in \Omega$, $|(\alpha) M_{(\Delta)} \setminus \Gamma|
\leq n$. For, otherwise $E_\leq \preq M$, by
Proposition~\ref{unbounded orbits}, and $E_{(A)} \prec E_\leq$, by
Theorem~\ref{fear}(ii).

Let $m > 1$ be the largest integer such that for all finite
$\Gamma, \Delta \subseteq \Omega$, there exists $\alpha \in
\Omega$ such that $|(\alpha) M_{(\Delta)} \setminus \Gamma| \geq
m$. Then there exist finite sets $\Gamma', \Delta' \subseteq
\Omega$ such that for all $\alpha \in \Omega$, $|(\alpha)
M_{(\Delta')} \setminus \Gamma'| \leq m$. Now, $M' =
M_{(\Delta')}$ has the property that for all finite $\Gamma,
\Delta \subseteq \Omega$ there exists $\alpha \in \Omega$ such
that $(\alpha)M'_{(\Delta)} \not\subseteq \Gamma \cup \{\alpha\}$,
and for all finite $\Gamma, \Delta \subseteq \Omega$, with
$\Gamma' \subseteq \Gamma$, the maximum of the cardinalities of
the sets $(\alpha) M'_{(\Delta)} \setminus \Gamma$ ($\alpha \in
\Omega)$ is $m$. From now on we will be working with $M'$ in place
of $M$.

A consequence of the above considerations is that if for some
$\alpha \in \Omega$ and finite $\Delta, \Gamma \subseteq \Omega$,
with $\Gamma' \subseteq \Gamma$, we have $|(\alpha) M'_{(\Delta)}
\setminus \Gamma| = m$, then $(\alpha) M'_{(\Delta)} \setminus
\Gamma = (\alpha) M' \setminus \Gamma'$, since $(\alpha) M'
\setminus \Gamma'$ cannot have cardinality larger than $m$. Thus
\begin{xlist}\item\label{unif bound2}
If $\Gamma, \Delta \subseteq \Omega$ are finite subsets,
with $\Gamma' \subseteq \Gamma$, and $\alpha \in \Omega$ has the
property that $|(\alpha) M'_{(\Delta)}\setminus\Gamma| = m$, then
for every $g \in M'$ that embeds $(\alpha) M'_{(\Delta)} \setminus
\Gamma$ in $\Omega \setminus \Gamma'$, we have $((\alpha)
M'_{(\Delta)} \setminus \Gamma)g = (\alpha) M'_{(\Delta)}
\setminus \Gamma$.
\end{xlist}
(This is because for such an element $g$, $((\alpha) M'_{(\Delta)}
\setminus \Gamma)g \subseteq (\alpha) M' \setminus \Gamma'$, and
the latter set is equal to $(\alpha) M'_{(\Delta)} \setminus
\Gamma$.)

We now construct recursively, for each $j \geq 0$, elements
$\alpha_j, \beta_j, \gamma_j \in \Omega$ and a finite subset $K_j
\subseteq M'$, consisting of a.e.\ injective maps, indexed
\begin{xlist}\item
$K_j = \{g(k_0, k_1, \ldots, k_{j-1}) : (k_0, k_1, \ldots,
k_{j-1}) \in \{0, 1\}^j\}.$
\end{xlist}
Let the $1$-element set $K_0 = \{g()\}$ consist of the identity
map $1 \in M'$. Assume inductively that for some $j \geq 0$ the
elements $\alpha_0, \dots, \alpha_{j-1}, \beta_0, \dots,
\beta_{j-1}, \gamma_0, \dots, \gamma_{j-1}$ and the sets $K_0,
\dots, K_j$ have been defined. Let $\Delta_j = \{\alpha_0, \ldots,
\alpha_{j-1}\} \cup \{1, \ldots, j-1\}$, and let $\Gamma_j
\subseteq \Omega$ be a finite set, containing $\Gamma' \cup
\{\beta_0, \ldots, \beta_{j-1}, \gamma_0, \ldots, \gamma_{j-1}\}$,
such that all elements of $K_0 \cup \dots \cup K_j$ embed $\Omega
\setminus \Gamma_j$ in $\Omega \setminus \Gamma'$. (Since $K_0
\cup \dots \cup K_j$ is finite and consists of elements that are
a.e.\ injective, we can find a finite set $\Gamma_j$, containing
$\Gamma' \cup \{\beta_0, \ldots, \beta_{j-1}, \gamma_0, \ldots,
\gamma_{j-1}\}$, such that all elements of $K_0 \cup \dots \cup
K_j$ are injective on $\Omega \setminus \Gamma_j$. We can then
enlarge this $\Gamma_j$ to include the finitely many elements that
are mapped to $\Gamma'$ by $K_0 \cup \dots \cup K_j$.) Now, let us
choose $\alpha_j \in \Omega$ so that $|(\alpha_j)M'_{(\Delta_j)}
\setminus \Gamma_j| = m$, and let $\beta_j, \gamma_j \in
(\alpha_j)M'_{(\Delta_j)} \setminus \Gamma_j$ be two distinct
elements. We note that $\beta_j$ and $\gamma_j$ are distinct from
$\beta_0, \ldots, \beta_{j-1}, \gamma_0, \ldots, \gamma_{j-1}$,
since the latter are elements of $\Gamma_j$.

Next, let us construct the elements $g(k_0, \ldots, k_{j-1}, 0),
g(k_0, \ldots, k_{j-1}, 1) \in K_{j+1}$. Writing $g' = g(k_0, k_1,
\ldots, k_{j-1}) \in K_j$, we define $g(k_0, \ldots, k_{j-1}, 0) =
hg'$ and $g(k_0, \ldots, k_{j-1}, 1) = fg'$, where $h, f \in
M'_{(\Delta_j)}$ are chosen so that $(\alpha_j)hg' = \beta_j$ and
$(\alpha_j)fg' = \gamma_j$. It is possible to find such $h$ and
$f$, by~(\ref{unif bound2}), using the fact that $g'$ embeds
$\Omega \setminus \Gamma_j$ in $\Omega \setminus \Gamma'$. By the
hypothesis that a.e.\ injective maps form a dense subset in $M$
(and hence in $M'$), we may further assume that $f$ and $g$ are
a.e.\ injective. As before, the images of $\alpha_0, \ldots,
\alpha_{j-1}$ and $1, \ldots, j-1$ under $hg'$ and $fg'$ will be
the same as their images under $g'$, and these two maps will be
a.e.\ injective (as composites of a.e.\ injective maps).

Given an infinite string $(k_i)_{i \in \omega} \in \{0,
1\}^\omega$, the elements of the sequence $g()$, $g(k_0)$, $\ldots
\,$, $g(k_0, \ldots, k_{i-1})$, $\ldots \,$ agree on successive
sets $\{1, \ldots, j-1\}$, and hence the sequence converges in
$E$. Since $M$ is closed, the limit belongs to $M$ (in fact,
$M'$). Considering our definition of the elements $\alpha_i$,
$\beta_i$, $\gamma_i$, the submonoid $M$ satisfies the hypotheses
of Lemma~\ref{tree3}, and hence $E_{(A)} \preq M$.
\end{proof}

In the previous lemma we were concerned with monoids under which
elements of $\Omega$, for the most part, had disjoint forward
orbits. Next, we give a similar argument, but adjusted for monoids
under which forward orbits can fall together.

\begin{lemma}
Let $M \subseteq E$ be a submonoid that is closed in the function
topology and has dense a.e.\ injective maps. If for every finite
$\Sigma \subseteq \Omega$, $M_{(\Sigma)} \neq \{1\}$, then
$E(\rho_\gamma) \preq M$ for some $\gamma \in \Omega$.
\end{lemma}

\begin{proof}
We may assume that there exist finite sets $\Gamma, \Delta
\subseteq \Omega$ such that for all $\alpha \in \Omega$, $(\alpha)
M_{(\Delta)} \subseteq \Gamma \cup \{\alpha\}$. For, otherwise,
the previous lemma implies that $E_{(A)} \preq M$, where $A$ is a
partition of $\Omega$ into $2$-element sets (and $E(\rho_\gamma)
\preq E_{(A)}$ for all $\gamma \in \Omega$, by
Lemmas~\ref{preord1}, \ref{preord2}, and~\ref{preord3}). Now, for
every finite $\Sigma \subseteq \Omega$, we have that
$(M_{(\Delta)})_{(\Sigma)} = M_{(\Delta \cup \Sigma)} \neq \{1\}$.
Also, $M_{(\Delta)}$ is closed in the function topology, and a.e.\
injective maps form a dense subset in $M_{(\Delta)}$. Hence, we
may replace $M$ with $M_{(\Delta)}$ and thus assume that
\begin{xlist}\item\label{fall together}
For every $\alpha \in \Omega$, $(\alpha) M \subseteq \Gamma \cup
\{\alpha\}$.
\end{xlist}
We now construct recursively for each $j \geq 0$ an element
$\alpha_j \in \Omega$ and a finite subset $K_j \subseteq M$,
consisting of a.e.\ injective maps, indexed
\begin{xlist}\item
$K_j = \{g(k_0, k_1, \ldots, k_{j-1}) : (k_0, k_1, \ldots,
k_{j-1}) \in \{0, 1\}^j\}.$
\end{xlist}
Let the 1-element set $K_0 = \{g()\}$ consist of the identity map
$1 \in M$. Assume inductively that for some $j \geq 0$ the
elements $\alpha_0, \dots, \alpha_{j-1}$ and the sets $K_0, \dots,
K_j$ have been defined. Since these sets are finite and consist of
a.e.\ injective maps, there is a finite set $\Delta_j \subseteq
\Omega$, such that all members of $K_0 \cup \dots \cup K_j$ are
injective on $\Omega \setminus \Delta_j$. We note that every $h
\in K_0 \cup \dots \cup K_j$, since it is injective on $\Omega
\setminus \Delta_j$, must act as the identity on all but finitely
many elements of $\Omega \setminus \Delta_j$, by~(\ref{fall
together}). Let $\Gamma_j \subseteq \Omega$ be the finite set
consisting of $\Delta_j \cup \Gamma \cup \{\alpha_0, \ldots,
\alpha_{j-1}\} \cup \{1, \ldots, j-1\}$ and those elements of
$\Omega \setminus \Delta_j$ that are taken to $\Gamma$ by members
of $K_0 \cup \dots \cup K_j$. We then take $\alpha_j \in \Omega
\setminus \Gamma_j$ to be any element such that
$|(\alpha_j)M_{(\Gamma_j)}| > 1$.

Next, let us construct the elements $g(k_0, \ldots, k_{j-1}, 0),
g(k_0, \ldots, k_{j-1}, 1) \in K_{j+1}$. Writing $g' = g(k_0, k_1,
\ldots, k_{j-1}) \in K_j$, we define $g(k_0, \ldots, k_{j-1}, 0) =
g'$, noting that $(\alpha_j)g' = \alpha_j$, by definition of
$\Gamma_j$. Also, let us define $g(k_0, \ldots, k_{j-1}, 1) =
fg'$, where $f \in M_{(\Gamma_j)}$ is chosen so that
$(\alpha_j)fg' \in \Gamma$ and $f$ is a.e.\ injective. (It is
possible to find an a.e.\ injective $f \in M_{(\Gamma_j)}$ such
that $(\alpha_j)f \in \Gamma$, by our hypotheses that
$|(\alpha_j)M_{(\Gamma_j)}| > 1$, that for every $\alpha \in
\Omega$, $(\alpha) M \subseteq \Gamma \cup \{\alpha\}$, and that
a.e.\ maps form a dense subset in $M$. Then $(\alpha_j)fg' \in
\Gamma$, by~(\ref{fall together}).) We note that the images of
$\alpha_0, \ldots, \alpha_{j-1}$ and $1, \ldots, j-1$ under $fg'$
will be the same as their images under $g'$, and that $fg'$ will
be a.e.\ injective.

Given an infinite string $(k_i)_{i \in \omega} \in \{0,
1\}^\omega$, the elements of the sequence $g()$, $g(k_0)$, $\ldots
\,$, $g(k_0, \ldots, k_{i-1})$, $\ldots \,$ agree on successive
sets $\{1, \ldots, j-1\}$, and hence the sequence converges in
$E$. Since $M$ is closed, the limit belongs to $M$. Now, let us
fix some $\gamma \in \Gamma$. Also, let $h_1 \in E$ be a map that
takes $\Omega$ bijectively to $\{\alpha_i\}_{i \in \omega}$, let
$h_2 \in E$ be the map that takes $\Gamma$ to $\gamma$ and fixes
all other elements of $\Omega$, and let $h_3 \in E$ be a right
inverse of $h_1$ that fixes $\gamma$. Then $E(\rho_\gamma)
\subseteq h_1 M h_2 h_3$.
\end{proof}

\begin{lemma}
Let $M \subseteq E$ be a submonoid that has large stabilizers.
Assume that there exist finite sets $\Gamma, \Delta \subseteq
\Omega$ such that for all $\alpha \in \Omega$, $(\alpha)
M_{(\Delta)} \subseteq \Gamma \cup \{\alpha\}$. Then $M \preq
E(\rho_\gamma)$ $($where $\gamma \in \Omega$ is any fixed
element$)$.
\end{lemma}

\begin{proof}
Since $M$ has large stabilizers, we may replace $M$ with
$M_{(\Delta)}$ and assume that for all $\alpha \in \Omega$,
$(\alpha) M \subseteq \Gamma \cup \{\alpha\}$. Let us define a
preorder $\rho$ on $\Omega$ by $(\alpha, \beta) \in \rho$ if and
only if $\beta \in (\alpha)M$. Then $M \subseteq E(\rho)$, and
$\rho$ is of type 3b. Hence, by Lemma~\ref{preord2}, $M \preq
E(\rho_\gamma)$.
\end{proof}

The following proposition summarizes the previous four lemmas.

\begin{proposition}\label{large stab 3}
Let $M \subseteq E$ be a submonoid that is closed in the function
topology, and has large stabilizers and dense a.e.\ injective
maps. Assume that there exists a finite set $\Sigma \subseteq
\Omega$ and a positive integer $n$ such that the cardinalities of
all forward orbits of $M_{(\Sigma)}$ are bounded by $n$, but for
every finite $\Sigma \subseteq \Omega$, $M_{(\Sigma)} \neq \{1\}$.
Let $A$ be a partition of $\Omega$ into $2$-element sets and
$\gamma \in \Omega$. Then
\begin{enumerate}
\item[\rm{(i)}] $M \approx E_{(A)}$ if and only if for
all finite $\Gamma, \Delta \subseteq \Omega$ there exists $\alpha
\in \Omega$ such that $(\alpha)M_{(\Delta)} \not\subseteq \Gamma
\cup \{\alpha\}$, and

\item[\rm{(ii)}] $M \approx E(\rho_\gamma)$ if and only if there
exist finite sets $\Gamma, \Delta \subseteq \Omega$ such that for
all $\alpha \in \Omega$, $(\alpha) M_{(\Delta)} \subseteq \Gamma
\cup \{\alpha\}$.
\end{enumerate}
\end{proposition}

\section{Main results}

Combining the results of the previous section, we obtain our
desired classification.

\begin{theorem}\label{classification}
Let $M \subseteq \E(\Omega)$ be a submonoid that is closed in the
function topology, and has large stabilizers and dense a.e.\
injective maps. Then $M$ falls into exactly one of five possible
equivalence classes with respect to $\approx$, depending on which
of the following conditions $M$ satisfies:
\begin{enumerate}
\item[\rm{(i)}] $M_{(\Sigma)}$ has an infinite forward orbit
for every finite $\Sigma \subseteq \Omega$.

\item[\rm{(ii)}] There exists a finite set
$\Sigma \subseteq \Omega$ such that all forward orbits of
$M_{(\Sigma)}$ are finite, but for every finite $\Sigma \subseteq
\Omega$ the cardinalities of the forward orbits of $M_{(\Sigma)}$
have no common finite bound.

\item[\rm{(iii)}] There exists a finite set
$\Sigma \subseteq \Omega$ and a positive integer $n$ such that the
cardinalities of all forward orbits of $M_{(\Sigma)}$ are bounded
by $n$, but for all finite $\Gamma, \Delta \subseteq \Omega$ there
exists $\alpha \in \Omega$ such that $(\alpha)M_{(\Delta)}
\not\subseteq \Gamma \cup \{\alpha\}$.

\item[\rm{(iv)}] There exist finite sets
$\Gamma, \Delta \subseteq \Omega$ such that for all $\alpha \in
\Omega$, $(\alpha) M_{(\Delta)} \subseteq \Gamma \cup \{\alpha\}$,
but for every finite $\Sigma \subseteq \Omega$, $M_{(\Sigma)} \neq
\{1\}$.

\item[\rm{(v)}] There exists a finite set
$\Sigma \subseteq \Omega$ such that $M_{(\Sigma)} = \{1\}$.
\end{enumerate}
\end{theorem}

\begin{proof}
Every submonoid of $E = \E(\Omega)$ clearly satisfies exactly one
of the above five conditions. By Corollary~\ref{large stab 1},
Corollary~\ref{large stab 2}, Proposition~\ref{large stab 3}, and
Theorem~\ref{fear}, the first four of these conditions describe
disjoint $\approx$-equivalence classes, when considering
submonoids that are closed in the function topology, and have
large stabilizers and dense a.e.\ injective maps.

Now, let $M, M' \subseteq E$ be submonoids that satisfy (v). Then
we can find finite sets $\Sigma, \Gamma \subseteq \Omega$ such
that $M_{(\Sigma)} = \{1\}$ and $M'_{(\Gamma)} = \{1\}$. If $M$
and $M'$ have large stabilizers, this implies that $M \approx
\{1\} \approx M'$. Finally, for submonoids of $E$ that satisfy the
hypotheses of the theorem, those that are described by condition
(v) are countable (since they are $\approx \{1\}$) and hence are
$\prec$ those submonoids that satisfy one of the other four
conditions.
\end{proof}

One can recover Theorem~\ref{partition classification} from this
by noting that when restricting to submonoids of the form
$E(\rho)$, the above five cases exactly correspond to the five
types of preorders on $\Omega$ identified in
Definition~\ref{types}.

If we instead restrict our attention to submonoids $M \subseteq E$
consisting of injective maps, then case (iv) of the above theorem
cannot occur, and case (iii) can be stated more simply.

\begin{corollary}\label{inj classification}
Let $M \subseteq \E(\Omega)$ be a submonoid that has large
stabilizers, is closed in the function topology, and consists of
injective maps. Then $M$ falls into exactly one of four possible
equivalence classes with respect to $\approx$, depending on which
of the following conditions $M$ satisfies:
\begin{enumerate}
\item[\rm{(i)}] $M_{(\Sigma)}$ has an infinite forward orbit for
every finite $\Sigma \subseteq \Omega$.

\item[\rm{(ii)}] There exists
a finite set $\Sigma \subseteq \Omega$ such that all forward
orbits of $M_{(\Sigma)}$ are finite, but for every finite $\Sigma
\subseteq \Omega$ the cardinalities of the forward orbits of
$M_{(\Sigma)}$ have no common finite bound.

\item[\rm{(iii)}] There exists
a finite set $\Sigma \subseteq \Omega$ and a positive integer $n$
such that the cardinalities of all forward orbits of
$M_{(\Sigma)}$ are bounded by $n$, but for every finite $\Sigma
\subseteq \Omega$, $M_{(\Sigma)} \neq \{1\}$.

\item[\rm{(iv)}] There exists a finite set $\Sigma \subseteq
\Omega$ such that $M_{(\Sigma)} = \{1\}$.
\end{enumerate}
\end{corollary}

Let $E^{\,\rm{inj}} \subseteq E$ be the submonoid consisting of
injective maps, and write $\approx_{\rm{inj}}$ to denote
$\approx_{\aleph_0, E^{\,\rm{inj}}}$. One might wonder whether the
above four $\approx$-classes are also
$\approx_{\rm{inj}}$-classes. They are not. For example, let
$E_\geq \subseteq E$ denote the submonoid of maps that are
increasing with respect to the usual ordering of $\omega = \Omega$
(i.e., maps $f \in E$ such that for all $\alpha \in \Omega$,
$(\alpha)f \geq \alpha$), and set $E^{\,\rm{inj}}_\geq =
E^{\,\rm{inj}} \cap E_\geq$. Then $E^{\,\rm{inj}}$ and
$E^{\,\rm{inj}}_\geq$ are both closed in the function topology (in
$E$), since limits of sequences of injective (respectively,
increasing) maps are injective (respectively, increasing). Also,
these two submonoids have large stabilizers. For, let $\Sigma
\subseteq \Omega$ be any finite set. Upon enlarging $\Sigma$, if
necessary, we may assume that $\Sigma = \{0, 1, \dots, n-1\}$ for
some $n \in \omega$. Let $f \in E$ be defined by $(i)f = i+n$ for
all $i \in \omega$, and let $g \in E$ be defined by $(i)g = i-n$
for $i \geq n$ (and arbitrarily on $\Sigma$). Then $E^{\,\rm{inj}}
\subseteq fE^{\,\rm{inj}}_{(\Sigma)}g$ and $E^{\,\rm{inj}}_\geq
\subseteq f(E^{\,\rm{inj}}_\geq)_{(\Sigma)} g$. (For, given $h \in
E^{\,\rm{inj}}$, respectively $E^{\,\rm{inj}}_\geq$, we can find
$\bar{h} \in E^{\,\rm{inj}}_{(\Sigma)}$, respectively
$(E^{\,\rm{inj}}_\geq)_{(\Sigma)}$, such that $(i + n)\bar{h} =
(i)h + n$ for all $i \in \omega$. Then $h = f\bar{h}g$.) This
shows that $E^{\,\rm{inj}}$ and $E^{\,\rm{inj}}_\geq$ satisfy the
hypotheses of the above corollary. They both clearly satisfy
condition (i), and hence $E^{\,\rm{inj}} \approx
E^{\,\rm{inj}}_\geq$. On the other hand, $E^{\,\rm{inj}}
\not\approx_{\rm{inj}} E^{\,\rm{inj}}_\geq$. For, suppose that
$E^{\,\rm{inj}} = \langle E^{\,\rm{inj}}_\geq \cup U \rangle$ for
some finite $U \subseteq E^{\,\rm{inj}}$. Then every element of
$E^{\,\rm{inj}}$ can be written as a word $f_0 f_1 \dots f_m$ in
elements of $E^{\,\rm{inj}}_\geq \cup U$. Let $f = f_0 f_1 \dots
f_{i-1}f_i f_{i+1} \dots f_m$ be such a word, and suppose that
$f_i \in E^{\,\rm{inj}}_\geq \setminus \{1\}$. As an increasing
injective non-identity element, $f_i$ is necessarily not
surjective (e.g., if $j \in \omega$ is the least element that is
not fixed by $f_i$, then $j \notin (\Omega)f_i$, since only
elements that are $\leq j$ can be mapped to $j$ by $f_i$). Since
$f_{i+1} \dots f_m$ is injective, this implies that $f$ cannot be
surjective either. Thus, if $f \in \langle E^{\,\rm{inj}}_\geq
\cup U \rangle$ is surjective, then $f \in \langle U \rangle$,
which is absurd, since $U$ is finite.

We can further show that $E^{\,\rm{inj}}$ and
$E^{\,\rm{inj}}_\geq$ satisfy the analog of the large stabilizer
condition in $E^{\,\rm{inj}}$ (i.e., for each finite subset
$\Sigma \subseteq \Omega$, $E^{\,\rm{inj}}_{(\Sigma)}
\approx_{\rm{inj}} E^{\,\rm{inj}}$ and
$(E^{\,\rm{inj}}_\geq)_{(\Sigma)} \approx_{\rm{inj}}
E^{\,\rm{inj}}_\geq$), demonstrating that the obvious analog of
the previous corollary for $E^{\,\rm{inj}}$ in place of $E$ is not
true. Let $\Sigma \subseteq \Omega$ be a finite subset. For each
injective map $f : \Sigma \rightarrow \Omega$, let us pick a
permutation $g_f \in E^{\,\rm{inj}}$ that agrees with $f$ on
$\Sigma$, and let $U$ be the (countable) set consisting of these
$g_f$. Now, let $f \in E^{\,\rm{inj}}$ be any map. Then $fg_f^{-1}
\in E^{\,\rm{inj}}_{(\Sigma)}$, which implies that $f \in \langle
E^{\,\rm{inj}}_{(\Sigma)} \cup U \rangle$. Now, as a countable set
of permutations, $U$ can be embedded in a subgroup of
$\Sym(\Omega)$ generated by two permutations, by
\cite[Theorem~5.7]{Galvin}, and hence in a submonoid of
$\Sym(\Omega) \subseteq E^{\,\rm{inj}}$ generated by four
elements. Therefore $E^{\,\rm{inj}} \approx_{\rm{inj}}
E^{\,\rm{inj}}_{(\Sigma)}$. To show that $E^{\,\rm{inj}}_\geq
\approx_{\rm{inj}} (E^{\,\rm{inj}}_\geq)_{(\Sigma)}$, we will
employ a similar method, though now we will assume, for
simplicity, that $\Sigma = \{0, 1, \dots, n-1\}$ for some $n \in
\Omega$, and require a more specific definition of the elements
$g_f$. For each injective $f : \Sigma \rightarrow \Omega$, let us
define a map $g_f \in E^{\,\rm{inj}}$ to agree with $f$ on
$\Sigma$, fix all elements not in $\Sigma \cup (\Sigma)f$, and act
on $(\Sigma)f$ in any way that turns $g_f$ into a permutation of
$\Omega$. As before, let $U$ be the set consisting of these
elements $g_f$. Now, let $f \in E^{\,\rm{inj}}_\geq$ be any map,
and let $h \in (E^{\,\rm{inj}}_\geq)_{(\Sigma)}$ be such that $h$
agrees with $f$ on $\Omega \setminus \Sigma$. We note that for any
$\alpha \in \Omega \setminus \Sigma$, $(\alpha)h \notin \Sigma
\cup (\Sigma)f$, since $f$ is injective and increasing. Hence $f =
hg_f$, and therefore $f \in \langle
(E^{\,\rm{inj}}_\geq)_{(\Sigma)} \cup U \rangle$. As before, this
implies that $E^{\,\rm{inj}}_\geq \approx_{\rm{inj}}
(E^{\,\rm{inj}}_\geq)_{(\Sigma)}$.

\section{Groups}\label{groups section}

We recall that the group $S = \Sym(\Omega)$ inherits from $E =
\E(\Omega)$ the function topology but is not closed in $E$ in this
topology. For instance, using cycle notation for permutations of
$\Omega = \omega$, we see that the sequence $(0,1), (0,1,2),
\ldots, (0, \ldots, n),\ldots$ converges to the map $i \mapsto
i+1$, which is not surjective. Thus we need a different notation
for the closure of a set of permutations of $\Omega$ in $S$; given
a subset $U \subseteq S$, let $\mathrm{cl}_S(U) = \mathrm{cl}_E
(U) \cap S$. It is easy to see that given a subset $U \subseteq
S$, $\mathrm{cl}_E (\mathrm{cl}_S (U)) = \mathrm{cl}_E (U)$.

Let $G \subseteq S$ be a subgroup, and let $\Sigma \subseteq
\Omega$ be finite. Then the forward orbits of $\mathrm{cl}_E
(G)_{(\Sigma)} = \mathrm{cl}_E (G_{(\Sigma)})$ coincide, by
Lemma~\ref{function top}, with the forward orbits of
$G_{(\Sigma)}$, which are simply the (group-theoretic) orbits of
$G_{(\Sigma)}$ in $\Omega$. In particular, by the above remark,
the forward orbits of $\mathrm{cl}_E (G)_{(\Sigma)} =
\mathrm{cl}_E (\mathrm{cl}_S (G))_{(\Sigma)}$ coincide with the
orbits of $\mathrm{cl}_S (G)_{(\Sigma)}$.

We are now ready to prove that while $\mathrm{cl}_S (G)$ and
$\mathrm{cl}_E (G)$ may be different for a given group $G$, they
are $\approx$-equivalent.

\begin{proposition}\label{closures}
Let $G$ be a subgroup of $S \subseteq E$. Then $\mathrm{cl}_S (G)
\approx \mathrm{cl}_E (G)$.
\end{proposition}

\begin{proof}
Let $A, B, C,$ and $D$ be partitions of $\Omega$ such that $A$
consists of only one set, $B$ consists of finite sets such that
there is no common finite upper bound on their cardinalities, $C$
consists of 2-element sets, and $D$ consists of 1-element sets. By
the main results of~\cite{B&S}, $\mathrm{cl}_S (G)$ is $\approx_S$
to exactly one of $S = S_{(A)}$, $S_{(B)}$, $S_{(C)}$, or $\{1\} =
S_{(D)}$ (see Theorem~\ref{monoids&groups} for the notation
$\approx_S$), and by Corollary~\ref{inj classification},
$\mathrm{cl}_E (G)$ is $\approx$ to exactly one of $E = E_{(A)}$,
$E_{(B)}$, $E_{(C)}$, or $\{1\} = E_{(D)}$. Moreover, by the above
remarks about (forward) orbits, $\mathrm{cl}_S (G) \approx_S
S_{(X)}$ for some $X \in \{A, B, C, D\}$ if and only if
$\mathrm{cl}_E (G) \approx E_{(X)}$. But, $S_{(X)} \approx
E_{(X)}$, by Proposition~\ref{partition monoids}, and
$\mathrm{cl}_S (G) \approx_S S_{(X)}$ if and only if
$\mathrm{cl}_S (G) \approx S_{(X)}$, by
Theorem~\ref{monoids&groups}. Hence $\mathrm{cl}_S (G) \approx
\mathrm{cl}_E (G)$.
\end{proof}

We have referred to the main results of Bergman and Shelah
in~\cite{B&S} while proving our results. However, it is
interesting to note that the Bergman-Shelah theorems can be
recovered from them.

\begin{theorem}[{\rm Bergman and Shelah}]
Let $G \subseteq \Sym(\Omega)$ be a subgroup that is closed in the
function topology on $S$. Then $G$ falls into exactly one of four
possible equivalence classes with respect to $\approx_S$ $($as
defined in Theorem~\ref{monoids&groups}$)$, depending on which of
the following conditions $G$ satisfies:
\begin{enumerate}
\item[\rm{(i)}] $G_{(\Sigma)}$ has an infinite orbit for
every finite $\Sigma \subseteq \Omega$.

\item[\rm{(ii)}] There exists a finite set $\Sigma \subseteq
\Omega$ such that all orbits of $G_{(\Sigma)}$ are finite, but for
every finite $\Sigma \subseteq \Omega$ the cardinalities of the
orbits of $G_{(\Sigma)}$ have no common finite bound.

\item[\rm{(iii)}] There exists
a finite set $\Sigma \subseteq \Omega$ and a positive integer $n$
such that the cardinalities of all orbits of $G_{(\Sigma)}$ are
bounded by $n$, but for every finite $\Sigma \subseteq \Omega$,
$G_{(\Sigma)} \neq \{1\}$.

\item[\rm{(iv)}] There exists a finite set $\Sigma \subseteq
\Omega$ such that $G_{(\Sigma)} = \{1\}$.
\end{enumerate}
\end{theorem}

\begin{proof}
Let $G_1, G_2 \subseteq \Sym(\Omega)$ be two subgroups that are
closed in the function topology. By the remarks at the beginning
of the section, $G_1$ and $G_2$ fall into the same one of the four
classes above if and only if $\mathrm{cl}_E (G_1)$ and
$\mathrm{cl}_E (G_2)$ fall into the same one of the four classes
in Corollary~\ref{inj classification} if and only if
$\mathrm{cl}_E (G_1) \approx \mathrm{cl}_E (G_2)$. By the previous
proposition, $\mathrm{cl}_E (G_1) \approx \mathrm{cl}_E (G_2)$ if
and only if $G_1 \approx G_2$, which, by
Theorem~\ref{monoids&groups}, occurs if and only if $G_1 \approx_S
G_2$.
\end{proof}

\section*{Acknowledgements}

The author is grateful to George Bergman, whose numerous comments
and suggestions have led to vast improvements in this note, and
also to the referee for mentioning related literature and
suggesting ways to simplify some of the arguments.

\vspace{.1in}

\noindent
Department of Mathematics \\
University of Southern California \\
Los Angeles, CA 90089 \\
USA \\

\noindent
Email: {\tt mesyan@usc.edu}
\end{document}